\documentclass{article}
\usepackage{layout}
\usepackage{amsthm}
\usepackage{amssymb}
\usepackage{latexsym}
\usepackage{euscript}
\usepackage{amsmath}
\usepackage{amstext}
\usepackage{amsgen}
\usepackage{amsbsy}
\usepackage{amsopn}
\usepackage{amsfonts}
\usepackage{dsfont}
\usepackage[frame,cmtip,curve,arrow,matrix,line,graph]{xy}
\usepackage{hyperref}
\usepackage{graphicx}

\DeclareMathAlphabet{\mathpzc}{OT1}{pzc}{m}{it}

\title{Witt groups of sheaves on topological spaces}
\author{Jon Woolf, Department of Mathematical Sciences,\\  Peach Street, University of Liverpool, L69 7ZL, UK.\\ email: Jonathan.Woolf@liv.ac.uk}
\date{August, 2007}

\newtheorem{theorem}{Theorem}[section]
\newtheorem{proposition}[theorem]{Proposition}
\newtheorem{corollary}[theorem]{Corollary}
\newtheorem{lemma}[theorem]{Lemma}

\theoremstyle{definition}
\newtheorem{definition}[theorem]{Definition}
\newtheorem{example}[theorem]{Example}
\newtheorem{examples}[theorem]{Examples}
\newtheorem{remark}[theorem]{Remark}

\newcommand{\ie}{i.e.\ }
\newcommand{\eg}{e.g.\ }

\newcommand{\zz}{\mathbb{Z}}
\newcommand{\qq}{\mathbb{Q}}

\newcommand{\cc}{\mathbb{C}}
\newcommand{\pp}{\mathbb{P}}

\newcommand{\im}{\mathrm{im}}

\newcommand{\cat}[1]{\mathbb{#1}}
\newcommand{\der}[1]{\mathbb{D}(#1)}
\newcommand{\perf}[1]{\mathbb{D}^p(#1)}
\newcommand{\perfc}[1]{\mathbb{D}^p_c(#1)}
\newcommand{\constr}[1]{\mathbb{D}^c(#1)}
\newcommand{\constrac}[1]{\mathbb{D}^{\mathbb{C}-c}(#1)}
\newcommand{\eqconstr}[2]{\mathbb{D}^{#1,c}(#2)}
\newcommand{\constrc}[1]{\mathbb{D}^c_c(#1)}
\newcommand{\ob}[1]{\mathcal{#1}}
\newcommand{\id}{id}
\newcommand{\shift}{\Sigma}
\newcommand{\mor}[2]{{\mathrm{Hom}}(#1,#2)}

\newcommand{\imor}[2]{\mathbf{Hom}(#1,#2)}

\newcommand{\bdy}{\partial}

\newcommand{\rshom}[2]{\mathbf{RHom}(#1,#2)}

\newcommand{\rdf}[1]{\mathbf{R}{#1}}
\newcommand{\struct}[1]{\ob{O}_{#1}}
\newcommand{\dual}{D}
\newcommand{\doubledual}{D^2}
\newcommand{\unit}{
\mathbf{1}}
\newcommand{\dualiser}{\mathcal{D}}
\newcommand{\funct}[2]{\mathrm{Funct}(#1,#2)}

\newcommand{\coh}[3]{H^{#1}(#2;#3)}
\newcommand{\ih}[2]{I\!H^{#1}(#2)}
\newcommand{\ic}[1]{\mathcal{I}\mathcal{C}(#1)}

\newcommand{\tensor}{\mathbf{\otimes}^L}

\newcommand{\commutingtriangle}[6]
{
\xymatrix{
{#1} \ar[rr]^{#4} \ar[dr]_{#5} && {#2} \\
& {#3} \ar[ur]_{#6}
}
}

\newcommand{\distinguishedtriangle}[6]
{
\xymatrix{
{#1} \ar[rr]^{#4}  && {#2} \ar[dl]_{#5}\\
& {#3} \ar@{-->}[ul]_{#6}
}
}

\newcommand{\sym}[2]{ \commutingtriangle{#1}{\dual
    #1}{\doubledual #1}{#2}{\Phi(#1)}{\dual #2 } }

\newcommand{\symch}[6]{
\xymatrix{ {#1} \ar[rr]^{#3} \ar[d]^{\Phi(#1)} &&
\dual {#1} \ar@{=}[d] \ar[rr]^{#4} &&
{#2} \ar[d]^{\partial {#3}} \ar[rr]^{#5} &&
\shift{#1} \ar[d]^{(-1)^r\Phi(\shift{#1})} \\
\doubledual {#1} \ar[rr]_{ \dual #3}  && 
\dual {#1} \ar[rr]_{\dual #5} &&
\dual {#2}  \ar[rr]_{{#4}} &&
\shift{#1} 
}}

\newcommand{\oct}[6]{ \xymatrix{ {#1} \ar@{-->}[dd] \ar[dr] && {#4}
    \ar@{-->}[ll] & {#1} \ar@{-->}[dd] && {#4} \ar@{-->}[ll]
    \ar@{-->}[dl] \\ & {#5} \ar@{-->}[dl] \ar[ur] & & & {#6} \ar[ul]
    \ar[dr] \\ {#2} \ar[rr] && {#3} \ar[ul] \ar[uu] & {#2} \ar[rr]
    \ar[ur] && {#3} \ar[uu] } }

\newcommand{\sh}[1]{\textrm{Sh}(#1)}
\newcommand{\shc}[1]{\textrm{Sh}^c(#1)}

\renewcommand{\hom}[3]{\textrm{Hom}_{#1}(#2,#3)}

\newcommand{\rx}{(R,K)}

\newcommand{\cochains}[2]{\Delta^{#1}(#2)}

\newcommand{\st}[1]{\textrm{st}({#1})}

\newcommand{\rxmod}{{\textrm{(R,K)-Mod}}}
\newcommand{\comrxmod}{{\textrm{Com(R,K)}}}
\newcommand{\rmod}{\textrm{R-Mod}}

\begin{document}

%\layout

\maketitle
\begin{abstract}
This paper investigates the Witt groups of triangulated categories of sheaves (of modules over a ring $R$ in which $2$ is invertible) equipped with Poincare-Verdier duality. We consider two main cases, that of perfect complexes of sheaves on locally compact Hausdorff spaces and that of cohomologically constructible complexes of sheaves on polyhedra. We show that the Witt groups of the latter form a generalised homology theory for polyhedra and continuous maps. Under certain restrictions on the ring $R$, we identify these constructible Witt groups of a finite simplicial complex with Ranicki's free symmetric L-groups. Witt spaces are the natural class of spaces for which the rational intersection homology groups have Poincare duality. When the ring $R$ is the rationals we identify the constructible Witt groups with the $4$-periodic colimit of the bordism groups of PL Witt spaces. This allows us to interpret L-classes of singular spaces as stable homology operations from the constructible Witt groups to rational homology.

MSC:		32S60, 19G99, 55U30, 57Q20.

Keywords:	Witt groups, Witt spaces, intersection cohomology, L-theory, topology of singular spaces.\end{abstract}

%\tableofcontents
\section{Introduction}

This paper investigates the Witt groups of triangulated categories of sheaves of $R$-modules equipped with Poincar\'e--Verdier duality. We will be most interested in the case when $R=\qq$, however the main results in \S 3 and \S 4 hold for any commutative regular Noetherian ring, of finite Krull dimension, in which $2$ is invertible. (By `regular' we mean that $R$ has finite global dimension and that every finitely generated module satisfies Auslander's condition.) We consider two main cases, that of perfect complexes of sheaves on locally compact Hausdorff spaces and that of cohomologically constructible complexes of sheaves on polyhedra. We show that the Witt groups of the latter, the constructible Witt groups, form a generalised homology theory for polyhedra and continuous maps. 

When every finitely generated $R$-module can be resolved by a finite complex of finitely generated free $R$-modules we identify the constructible Witt groups of a finite simplicial complex $K$ with Ranicki's free symmetric L-groups $H_*(K; \mathbb{L}^\cdot(R))$ \cite[Proposition 14.5]{ranicki}. 

When $R=\qq$ we show that every Witt space has a natural L-theory, or Witt, orientation and we identify the constructible Witt groups with the $4$-periodic colimit of the bordism groups of Witt spaces introduced in \cite{siegel}. This answers Problem 6 in \cite[\S IX]{bo}. It also allows us to interpret L-classes of singular spaces as stable homology operations from the constructible Witt groups to rational homology. Before giving further details we put these results into context.

\subsection*{Witt groups and L-theory}

 In his 1937 paper \cite{witt} Witt studied symmetric bilinear forms over a field $k$, in particular defining what is now known as the Witt group $W(k)$ --- the set of isometry classes of symmetric bilinear forms (equipped with direct sum) modulo the stable equivalence relation generated by those forms with a Lagrangian subspace. (Witt also showed that the tensor product gives $W(k)$ a natural ring structure, but we will ignore this for the present.) By analogy we can define the Witt group $W(R)$ of any commutative ring $R$ --- see \eg \cite{knus,mh}. This algebraic construction has been generalised to provide invariants in both algebraic geometry and in algebraic topology. 
 
In algebraic geometry, Knebusch defined the Witt group $W(S)$ of a scheme $S$ in \cite{knebusch} by considering symmetric bilinear forms on locally-free coherent sheaves (vector bundles) on $S$. In this context the classical Witt group $W(R)$ of a ring $R$ arises as $W(\mathrm{Spec}\, R)$. Knebusch's definitions can be used to define the Witt group of any exact category with duality. In a more recent development \cite{balmer} Balmer extended this to define the Witt groups of any triangulated category $\cat{T}$ with duality. To obtain a good theory he requires that $2$ be invertible \ie that the morphisms between any two objects are a $\zz[\frac{1}{2}]$-module not merely an Abelian group. Balmer's Witt groups are a collection $W^i(\cat{T})$ of Abelian groups indexed by $\zz$, but which turn out to be naturally $4$-periodic i.e. $W^i(\cat{T}) \cong W^{i+4}(\cat{T})$. In a series of papers Balmer and others, notably Gille
and Walter, have studied these groups for the derived
category $\mathbb{D}^\textrm{lf}(S)$ of locally-free coherent
sheaves on a scheme. Knebusch's Witt group $W(S)$ is isomorphic to Balmer's zero'th Witt group $W^0(\mathbb{D}^\textrm{lf}(S))$.  Much of this work is
summarised in \cite[\S 5]{handbook}, which also contains a
compendious bibliography. Of particular note is \cite{hornbostel} in
which it is shown that the Witt groups of the derived category of locally-free coherent
sheaves on a regular scheme are representable in both the stable and unstable
$\mathbb{A}^1$-homotopy categories. 

In algebraic topology, the development by Browder, Novikov, Sullivan and Wall of the surgery theory of
high-dimensional manifolds in the 1960s culminated in the introduction by Wall \cite{wall} of the surgery obstruction groups
$L_*(R)$. These L-groups are defined for any ring with involution R and are 4-periodic \ie $L_*(R)\cong L_{*+4}(R)$.
Mishchenko and Ranicki also defined symmetric L-groups $L^*(R)$, with $L^0(R)=W(R)$ the classical Witt
group $W(R)$. If $2$ is invertible in $R$ then $L_i(R)\cong L^i(R)$, and 
$$L_i(R) \cong W^i(D^\textrm{lf}(\textrm{Spec}\, R)).$$
The L-groups of a ring with involution can be interpreted as the invariants associated to a point in a wider topological theory. More precisely,  given a simplicial complex $K$ and ring $R$ with involution, Ranicki has shown how to associate to it certain algebraic objects, called $(R,K)$-modules. Using a combinatorial version of Poincar\'e duality he constructs a `weak chain duality' on the category of chain complexes of $(R,K)$-modules and defines  the symmetric L-groups of $K$ to be the `algebraic bordism groups' of `symmetric Poincar\'e complexes' in this category. Furthermore he constructs a symmetric L-theory spectrum $\mathbb{L}^\cdot(R)$ whose homotopy groups are the symmetric L-groups of $R$. This spectrum corresponds to a generalised homology theory whose homology groups $H_*(K;\mathbb{L}^\cdot(R))$ are the symmetric L-groups of  $K$. L-theory plays an important r\^ole in surgery theory and the classification of manifolds. The definitive account of this work is \cite{ranicki}. 

In short, in both algebraic geometry and algebraic topology one can define generalised homology theories (in the sense appropriate to each subject) for which the classical Witt group appears as the zero'th group of a point. It is important to realise that the dualities involved are rather different in these two cases; in algebraic geometry one only has to extend vector space (free $R$-module) duality to vector bundles (locally-free coherent sheaves), but in topology one requires some form of Poincar\'e duality. There is another way to extend Witt groups in algebraic geometry where, rather than considering vector bundles, one considers the derived category of coherent sheaves equipped with Serre duality --- see \cite[Example 5.2]{handbook}. 

This paper draws from both these theories in that we apply Balmer's techniques, which arose in algebraic geometry, to obtain a new description of symmetric L-theory for polyhedra. 

\subsection*{Survey of results}

In slightly more detail, the contents of the paper are as follows. We begin, in  \S\ref{categorical witt groups}, by surveying the basic definitions and properties of Balmer's Witt groups of a triangulated category with duality. We elucidate the connection between symmetric isomorphisms in a triangulated category with duality and Verdier dual pairings, which are the analogue of symmetric bilinear forms. This section also contains a new treatment of the appropriate functors between categories with duality, namely functors which are symmetrically self-dual. 

We apply this theory in \S\ref{perfect witt groups} to construct Witt groups of sheaves on topological spaces. Suppose $X$ is a locally compact, locally connected Hausdorff space which is countable at infinity and $R$ a commutative regular Noetherian ring of finite Krull dimension. Under these conditions there is a (contravariant) Poincar\'e--Verdier duality functor from the derived category of sheaves of $R$-modules on $X$ to itself (see \cite[Chapter 3]{dimca}). If we restrict to the triangulated subcategory of perfect complexes then it is an equivalence, and we have a triangulated category with duality in the sense of \cite{twg1}. With the further assumption that $2$ is invertible, we show that its Witt groups $W^p_*(X)$ form a homotopy-invariant functor which satisfies all the axioms of a generalised homology theory apart from possibly excision. 

Let $K$ be  a simplicial complex. Its realisation is naturally stratified (see \S\ref{comb and constr}), and we denote this stratified space by $K_S$.  We can restrict our attention from the perfect complexes to the triangulated subcategory of  complexes which are cohomologically constructible with respect to the stratification. This subcategory is also preserved by Poincar\'e--Verdier duality. Its Witt groups, which we dub the constructible Witt groups of $K$ and denote $W^c_*(K)$, form a generalised homology theory for simplicial complexes and simplicial maps. Simplicial approximation then allows us to obtain a generalised homology theory for compact polyhedra and continuous maps. The section ends with a brief discussion of  equivariant
generalisations and of related theories defined by altering the
constructibility condition.

In  \S\ref{L-theory} we relate the constructible Witt groups $W_*^c(K)$ of a finite simplicial complex $K$ to Ranicki's free symmetric L-groups $H_*(K; \mathbb{L}^\cdot(R))$ by exhibiting a natural transformation from the latter to the former. If every finitely generated $R$-module has a finite resolution by finitely generated free $R$-modules then a theorem of Walter's \cite[Theorem 5.3]{walter} shows that the natural transformation induces an isomorphism of point groups.  Hence we obtain isomorphic generalised homology theories for simplicial complexes.

Finally, \S\ref{rational theory} explains the geometric nature of the rational ($R=\qq$) theory.  We review Siegel's work \cite{siegel} on the bordism
groups of PL Witt spaces and construct a natural transformation from Witt bordism to
the constructible Witt groups which is, in sufficiently high dimensions, an
isomorphism. Phrased another way,  the constructible Witt groups are the $4$-periodic colimit of the Witt bordism groups. We use this geometric description, and
an adaptation of the construction of L-classes in \cite[\S 20]{milstaff}, to view L-classes as homology operations
from the bordism groups of Witt spaces, or, by the identification of
the previous section, from the constructible Witt groups, to
rational homology. 

\subsection*{Connections with other work}

The isomorphism between certain constructible Witt groups and free symmetric L-groups constructed in \S\ref{L-theory}  makes it apparent that this paper is closely related to Ranicki's work on L-theory.  Our sheaf-theoretic approach has the virtues that it is
technically simpler (at least for those familiar with sheaves and
derived categories) and that it directly connects L-theory with the large body
of work on intersection homology, self-dual complexes of sheaves and
characteristic classes for singular spaces. 

This is not the first attempt to give a sheaf-theoretic description of L-theory. One could loosely describe sections \ref{perfect witt groups} and \ref{L-theory} as a
triangulated version of Hutt's unpublished paper \cite{hutt}, in which he considers the symmetric and quadratic L-groups of the category of complexes of sheaves with Poincar\'e--Verdier duality. However, there are important differences. By working directly with complexes, rather than in the derived category, Hutt obtains distinct quadratic and symmetric L-groups. Using the triangulated approach we require the restriction that $2$ is invertible. (Although it is possible to define Witt groups for triangulated categories of sheaves even when $2$ is not invertible, there is no known proof of the long exact sequence of a pair in this case, and it seems doubtful that we would obtain a good theory, see \cite[\S 4]{handbook}.) This means that we obtain only one theory, since the quadratic and symmetric L-groups agree when $2$ is invertible. From the point of
view of topology, in which we are most interested in
$\mathbb{L}^\cdot(\zz)$ or $\mathbb{L}^\cdot(\zz[\pi_1(X)])$, the
restriction to rings in which $2$ is invertible is
perhaps unfortunate. However, in compensation for this restriction, the triangulated theory
is considerably simpler and less fiddly to define. The proofs are quite formal, often based on nothing more than well-known
properties of functors between triangulated categories of sheaves. A
case in point is the proof of excision for Witt groups of
constructible sheaves in \S\ref{excision}. We do not require any of the
machinery of micro-supports involved in Hutt's work, and thus we are
able to obtain a
generalised homology theory for polyhedra and continuous maps without requiring any assumptions of smoothness.

There is also a close relation, particularly in terms of technique, to
Cappell and Shaneson's work on self-dual complexes of sheaves, see
\cite{cs}, and to Youssin's more formal version \cite{youssin} for
triangulated categories, in which he similarly defines a cobordism relation on self-dual objects in a derived category with duality.  However, this
cobordism relation, see \cite[\S 2]{cs} for the case of sheaves and
more generally \cite[Definition 6.1]{youssin}, is stronger than the
relation of Witt-equivalence introduced by Balmer \cite[1.13]{twg1}
which we use. Another difference is that Cappell and Shaneson work
with a fixed stratification. For this reason one would not expect their cobordism groups to form a generalised homology theory. To obtain a topologically invariant
theory we take a limit over all stratifications compatible with a
given PL structure. Taking into
account these differences, and using \cite[Remark 3.25]{twg1}, we see that Cappell and Shaneson's cobordism group is a quotient of
the Witt group of the triangulated subcategory of complexes of sheaves which are  \emph{constructible with respect to
  a fixed stratification}. In particular, their cobordism
groups are always freely generated Abelian groups (see \cite[Theorem
  4.7]{cs} and \cite[Corollary 7.5]{youssin}) whereas ours can, and
frequently do, have torsion. The relationship between Balmer's Witt groups and Youssin's cobordism groups is discussed in more detail in \cite{chern-classes}. (A potential source of confusion is that Youssin calls his cobordism groups  Witt groups. However, it should be noted that these Witt groups are \emph{not} the same as Balmer's and, moreover, Youssin's cobordism group of the derived category of $R$-modules is \emph{not} in general the Witt group of the ring $R$.)

A large part of \cite{cs}, and also of \cite{bcs}, is concerned with
the computation of (intersection cohomology) signatures and, more
generally, L-classes, which they show are invariant under their
cobordism relation \cite[Proposition 5.2]{cs}. Since their cobordism
group is a quotient of our Witt group  (provided we use a
fixed stratification) it follows that
L-classes are also well defined on the Witt group. We give a direct geometric
construction in \S\ref{L-classes} showing that the L-classes 
arise as stable homology operations from the constructible Witt groups  to ordinary rational homology. See also the construction of L-classes for singular varieties in \cite{chern-classes}.

Further information on the connections between surgery and L-theory, self-dual sheaves and Witt spaces can be found in the excellent survey \cite{hw}. The reader is also referred to \cite{banagl}, particularly \S 4, in which a bordism group $\Omega_*^{SD}$ of self-dual sheaves is constructed. The bordism relation has both a geometric and an algebraic component, the latter of which is similar to the cobordism relation in \cite{cs}. There is a natural map from the bordism groups of Witt spaces to $\Omega_*^{SD}$ (through which the signature factorises) but it is not immediately clear how the two theories are related.

\subsection*{Acknowledgments}
I would like to thank Andrew Ranicki for suggesting the topic of this paper to me and for patiently reading preliminary drafts and correcting my many misapprehensions. I would also like to thank the referee and J\"org Sch\"urmann for their helpful comments and corrections.

\section{Balmer's Witt Groups}
\label{categorical witt groups}
\subsection{Definitions}
Let $\cat{A}$ be a triangulated category. We will assume that it satisfies the \emph{enriched octahedral axiom} discussed in \cite[Remarque 1.1.13]{bbd} and \cite{balmer}. As noted in these references, this is satisfied by all known examples of triangulated categories, in particular by derived categories. It also passes to triangulated subcategories and localisations. The shift functor will be denoted $\shift$ and exact triangles written  $\ob{A} \to \ob{B} \to\ob{C}\to\shift\ob{A}$. Throughout  we assume that $2$ is invertible \ie given $\alpha \in \mor{\ob{A}}{\ob{B}}$ there exists $\alpha'$ with $\alpha = 2 \alpha'$.

Following Balmer \cite[\S1.2]{twg1} we say a pair $(\dual, \Phi)$ is a \emph{$\delta$-duality}, where $\delta = \pm 1$, if 
\begin{enumerate}
\item $\dual$ is an additive functor $\dual :
\cat{A} \to \cat{A}^{op}$  such that for any exact triangle $\ob{A} \stackrel{\alpha}{\longrightarrow} \ob{B} \stackrel{\beta}{\longrightarrow}\ob{C}\stackrel{\gamma}{\longrightarrow}\shift\ob{A}$ the triangle 
$\dual\ob{C} \stackrel{\dual\beta}{\longrightarrow} \dual\ob{B} \stackrel{\dual \alpha}{\longrightarrow} \dual\ob{A}\stackrel{\delta\dual\gamma}{\longrightarrow}\shift\dual\ob{C}$
 is exact and
\item $\Phi$ is an isomorphism of functors $\id \stackrel{\Phi}{\longrightarrow}{\dual}^2$ with $\Phi_{\shift \ob{A}} = \shift \Phi_A$ and satisfying the coherence relation 
\begin{equation}
\label{coherence}
\dual \Phi_\ob{A} \circ \Phi_{\dual\ob{A}} = \id_{\dual\ob{A}}
\end{equation}
(from which it follows that $\Phi_{\dual\ob{A}} \circ \dual \Phi_\ob{A} = \id_{\dual^3\ob{A}}$). We will give another interpretation of this coherence relation in Example \ref{symmetry of phi}.
\end{enumerate}
%In the examples we consider in \S\ref{perfect witt groups} and \ref{rational theory} there will be a %natural $1$-duality \ie a triangulated endofunctor $\dual$ whose square is isomorphic to the identity %via $\Phi$. 

The shift functor in a triangulated category is additive but not triangulated; we must change an odd number of signs of morphisms to regain an exact triangle from the shift of an exact triangle. With this in mind, it is not difficult to see that $(\shift\dual, -\delta\Phi)$ is a $(-\delta)$-duality. We will call this the \emph{shifted duality}. Note that the $r^{th}$ shifted duality of a $1$-duality $(D,\Phi)$ is a $(-1)^r$-duality given by 
$$
(\shift^r\dual, (-1)^{r(r+1)/2}\Phi).
$$
(In fact it is easy to check that whenever $(\dual, \Phi)$ is a $\delta$-duality then so is $(\dual, -\Phi)$ so the sign $(-)^{r(r+1)/2}\Phi$ is purely conventional. Nevertheless, it turns out to be the more natural choice for reasons which will become apparent in Lemma \ref{boundary lemma} below.) These shifted dualities will give us the higher Witt groups. 

Fix a $1$-duality $(\dual,\Phi)$ on $\cat{A}$. A morphism $\ob{A} \stackrel{\alpha}{\longrightarrow} \dual\ob{A}$ is \emph{symmetric} if the diagram
$$
\sym{\ob{A}}{\alpha}
$$
commutes. More generally, we say a morphism is \emph{symmetric of dimension $r$} if it is symmetric for the $r^{th}$ shifted duality. If $\alpha$ is an isomorphism we say that $\ob{A}$ is \emph{symmetrically self-dual via} $\alpha$. Symmetric morphisms $\alpha$ and $\beta$ are said to be \emph{isometric} if there is a commutative diagram
\[
\xymatrix{
\ob{A} \ar[d]_{\eta} \ar[r]^{\alpha} & \dual\ob{A} \\
\ob{B} \ar[r]_{\beta} & \dual\ob{B} \ar[u]_{\dual\eta}
}
\]
in which $\eta$ is an isomorphism. 
\begin{lemma}[{Balmer \cite[1.6]{twg1}}]
\label{boundary lemma}
Let $\alpha: \ob{A} \to \dual\ob{A}$ be a symmetric morphism. Then for any choice $\ob{B}$ of cone  on $\alpha$ we can choose a symmetric morphism $\beta$ of dimension $1$ such that  
$$
%\symch{\ob{A}}{\ob{B}}{\alpha}{\beta}{\gamma}{r}\footnote{sort out!}
\xymatrix{
\ob{A} \ar[d]_{\phi(\ob{A})} \ar[rr]^{\alpha} && \dual \ob{A} \ar@{=}[d] \ar[rr]^{\alpha'}&& \ob{B}\ar[d]^{\beta} \ar[rr]^{\alpha''}&& \shift \ob{A}  \ar[d]^{\shift\phi(\ob{A})}\\
\doubledual \ob{A} \ar[rr]_{\dual\alpha}&& \dual \ob{A} \ar[rr]_{-\shift\dual\alpha''}&& \shift \dual\ob{B} \ar[rr]_{\shift\dual\alpha'}&& \shift \doubledual \ob{A}
}
$$
(whose rows are exact triangles) commutes. Furthermore, if $\ob{B}'$ and $\beta'$ are different choices for the cone on $\alpha$ and the completing symmetric morphism, then $\beta$ and $\beta'$ are isometric.
\end{lemma} 
Note that $\beta$ is always an isomorphism so that repeating the coning construction starting with $\beta$ gives zero.

We define the $r^{th}$ monoid of symmetric morphisms $Sym^r(\cat{A},\dual,\Phi)$ to be the set of isometry classes of symmetric morphisms of dimension $r$ equipped with the addition arising from direct sum. Taking the cone of a symmetric morphism defines a coboundary operator
$$
Sym^r(\cat{A},\dual,\Phi) \stackrel{d}{\longrightarrow} Sym^{r+1}(\cat{A},\dual,\Phi)
$$
with $d^2=0$.
\begin{definition}[{Cf.\ Balmer \cite[1.13]{twg1}}]
The {\em $r^{th}$ Witt group} $W^r(\cat{A},\dual,\Phi)$ is the quotient of the monoid $$\ker \left( d : Sym^r(\cat{A},\dual,\Phi) \to Sym^{r+1}(\cat{A},\dual,\Phi) \right)$$ by the submonoid $\im \left( d : Sym^{r-1}(\cat{A},\dual,\Phi) \to Sym^{r}(\cat{A},\dual,\Phi) \right)$. Usually we suppress $\dual$ and $\Phi$ and simply write $W^r(\cat{A})$. The quotient is a group, with $-\alpha$ representing the additive inverse of the class of $\alpha$. 
\end{definition}

The Witt groups of any triangulated category with duality are naturally $4$-periodic. Whenever $\alpha : \ob{A} \to \dual\ob{A}$ is a symmetric morphism for $(\dual,\Phi)$ then the shift $\shift\alpha : \shift\ob{A} \to \shift\dual\ob{A} = \shift^2\dual\shift\ob{A}$ is symmetric for $(\shift^2\dual,\Phi)$ which is the second shifted duality of $(\dual, -\Phi)$. This defines an isomorphism 
$$
W^r(\cat{A},\dual,\Phi) \cong W^{r+2}(\cat{A},\dual,-\Phi)
$$
which, when repeated, yields the 4-periodicity
$W^r(\cat{A},\dual,\Phi) \cong W^{r+4}(\cat{A},\dual,\Phi)$.
\subsection{Internal structures and symmetric forms}
\label{internal structures}

In linear algebra we are familiar with the correspondence between (symmetric) bilinear forms and (symmetric) maps from a vector space to its dual. A similar interpretation is possible for the symmetric morphisms defined above, provided we add extra structure to our triangulated category.

First, we require that $\cat{A}$ be a symmetric monoidal category. In other words we have an (additive) tensor product 
$\otimes : \cat{A} \times \cat{A} \longrightarrow \cat{A}$
and functorial isomorphisms $\sigma_{\ob{AB}}:A \otimes B \to B\otimes A$, and there is a unit $\unit\in\cat{A}$ with $\unit \otimes \ob{A} \cong \ob{A}$ for all $\ob{A}$. Second, there should be an internal hom functor
$$
\imor{-}{-} : \cat{A}^{op} \times \cat{A} \to \cat{A}
$$
which is compatible with the tensor product in that
$$
\imor{\ob{A}}{\imor{\ob{B}}{\ob{C}}} \cong \imor{\ob{A}\otimes\ob{B}}{\ob{C}}
$$
for all $\ob{A},\ob{B}$ and $\ob{C}$. This should be related to morphisms in $\cat{A}$ by a functor $\Gamma$ from $\cat{A}$ to Abelian groups with $\mor{-}{-} \cong \Gamma \circ \imor{-}{-}$.

Finally we require that the duality is \emph{internal} \ie it is represented by a dualising object $\dualiser$ with respect to the internal hom, so that $\dual\ob{A} = \imor{\ob{A}}{\dualiser}$. 

\begin{remark}
We can express the internal hom in terms of the tensor product and duality as
\begin{equation}
\label{internal hom via tensor prod}
\imor{\ob{A}}{\ob{B}} \cong \imor{\ob{A}}{\dual^2\ob{B}} \cong  \imor{\ob{A} \otimes \dual \ob{B}}{\dualiser} \cong\dual(\ob{A}\otimes \dual\ob{B}).
\end{equation}
An immediate consequence is that $\imor{\unit}{\ob{A}} \cong\dual(\unit\otimes \dual\ob{A}) \cong \dual^2\ob{A} \cong \ob{A}$. In particular the dualising object is isomorphic to the dual of the unit: $\dualiser \cong \imor{\unit}{\dualiser} \cong \dual \unit$. 
\end{remark}

In the presence of this extra structure we see that there is an isomorphism
$$
\Theta:\mor{\ob{A}}{\dual\ob{B}} \cong \mor{\ob{A}}{\imor{\ob{B}}{\dualiser}} \cong \mor{\ob{A}\otimes\ob{B}}{\dualiser}.
$$
 The following lemma expresses the key properties of this correspondence.
 \begin{lemma}
 The following diagrams commute:
 $$
 \xymatrix{
 \ob{A}\otimes\ob{B} \ar[dr]^{\Theta(\alpha)} \ar[d]_{\sigma_{\ob{AB}}} && \ob{A}\otimes\ob{B} \ar[dr]^{\Theta(\dual\alpha\circ\gamma\circ\beta)}\ar[dd]_{\alpha\otimes\beta}&& \dual \ob{A} \otimes \ob{A} \ar[dr]^{\Theta(\id_{\dual\ob{A}})} \ar[dd]_{\sigma_{\dualiser\ob{A}\ob{A}}}&\\
 \ob{B}\otimes\ob{A} \ar[d]_{\Phi \otimes \id} & \dualiser&& \dualiser&&\dualiser.\\
 \doubledual\ob{B}\otimes\ob{A} \ar[ur]_{\Theta(\dual\alpha)}&&\ob{A}'\otimes\ob{B}' \ar[ur]_{\Theta(\gamma)}&& \ob{A} \otimes \dual \ob{A} \ar[ur]_{\Theta(\Phi\ob{A})}
 }
 $$
 \end{lemma}
\begin{proof}Exercise!
\end{proof}
It follows that this correspondence takes symmetric morphisms to symmetric bilinear forms, that is $\beta \in \mor{\ob{A}\otimes\ob{A}}{\dualiser}$ with $\beta \circ \sigma_{\ob{AA}} = \beta$. Symmetric \emph{isomorphisms} in $\mor{\ob{A}}{\dual\ob{A}} $ correspond to \emph{non-degenerate} symmetric bilinear forms, which, in this context, are forms $\beta$ with the property that 
$$
\beta \circ (\gamma \otimes \delta) = 0 \ \textrm{for all}\  \delta \iff \gamma = 0.
$$
We will also say that a pairing $\ob{A}\otimes\ob{B} \to \dualiser$ which corresponds to an isomorphism $\ob{A} \to \dual\ob{B}$ is a \emph{Verdier dual pairing}. In particular, symmetric isomorphisms yield Verdier dual pairings.

\begin{remark}
In \S\ref{perfect witt groups}, when we look at Witt groups of sheaves on a topological space, the relevant triangulated categories will possess all of this additional structure. 

In linear algebra we are used to the situation in which the unit $\unit$ and dualising object $\dualiser$ are naturally isomorphic, so that a bilinear form is a map from the tensor product of an object with itself to the unit. Furthermore there are natural isomorphisms
$$
\dual(\ob{A}\otimes\ob{B}) \cong \dual\ob{A} \otimes \dual\ob{B}.
$$
so that (\ref{internal hom via tensor prod}) becomes the more familiar $\imor{\ob{A}}{\ob{B}} \cong \dual\ob{A} \otimes \ob{B}$. Neither of these further properties will hold in the examples we consider in \S\ref{perfect witt groups}.
\end{remark}

\subsection{Functors and duality}
\label{functors and duality}
What functorial properties do the Witt groups have? A general triangulated functor $F : \cat{A} \to \cat{B}$ between triangulated categories with $\delta$-dualities will not preserve symmetric morphisms, and so cannot be expected to induce a map of Witt groups. The functor should be `symmetrically self-dual' too. 

Note that the functor category $\funct{\cat{A}}{\cat{B}}$ whose objects are triangulated functors from $\cat{A}$ to $\cat{B}$ and whose morphisms are natural transformations inherits a $\delta$-duality
\begin{equation}\label{functor duality}
\dual_{\cat{A},\cat{B}} F := \dual_\cat{B} \circ F \circ \dual_\cat{A}
\end{equation}
where $\delta_{\cat{A},\cat{B}} = \delta_\cat{A}\delta_\cat{B}$. There is an associated natural transformation (of morphisms of functors)
$\Phi_{\cat{A},\cat{B}} : \id \to \doubledual_{\cat{A},\cat{B}}$
which applied to a functor $F$ is the natural transformation 
\begin{equation}
\label{functor phi}
\Phi_{\cat{A},\cat{B}} (F) (-) = \Phi_\cat{B}(F \circ \doubledual_\cat{A} (-) ) \circ F ( \Phi_\cat{A} (-)).
\end{equation}
Note that using the commutative square
$$
\xymatrix{
F \ar[r]^{F\Phi_\cat{A}} \ar[d]_{\Phi_\cat{B}(F)} & F\dual_\cat{A}^2
\ar[d]^{\Phi_\cat{B}(F\dual_\cat{A}^2)} \\
\dual_\cat{B}^2F\ar[r]_{\dual_\cat{B}^2F\Phi_\cat{A}} & \dual_\cat{B}^2F\dual_\cat{A}^2
}
$$
we also have $\Phi_{\cat{A},\cat{B}} (F) (-) = \dual_\cat{B}^2 F(\Phi_\cat{A}(-)) \circ \phi_\cat{B}(F(-))$.

\begin{remark}
Here we must be careful about the meaning of $\delta$-duality since $\funct{\cat{A}}{\cat{B}}$ is not naturally triangulated. (This unfortunate situation arises because the cone on a morphism is not functorial, or, put another way, there is not necessarily a functor which is the cone on a morphism of functors.) However, $\funct{\cat{A}}{\cat{B}}$ is additive, has a well-defined shift operator and we can identify `exact triangles' to be diagrams of functors $F \to G \to H \to \shift F$ such that
$$
F\ob{A} \to G\ob{A} \to H\ob{A} \to \shift F\ob{A}
$$
is an exact triangle in $\cat{B}$ for any $\ob{A} \in \cat{A}$.
\end{remark}

We define a \emph{symmetric morphism of functors} to be a morphism of functors, \ie a natural transformation, $\Psi:F \to \dual_{\cat{A},\cat{B}} F$ such that
$$
 \commutingtriangle{F}{\dual_{\cat{A},\cat{B}} F}{\doubledual_{\cat{A},\cat{B}} F}{\Psi}{\Phi_{\cat{A},\cat{B}}F}{\dual_{\cat{A},\cat{B}} \Psi }
$$
commutes in $\funct{\cat{A}}{\cat{B}}$. As before we define symmetric morphisms of functors of non-zero dimension using the shifted dualities on  $\funct{\cat{A}}{\cat{B}}$.
\begin{remark}
The definition of a symmetric isomorphism of functors is equivalent to Balmer's definition of a morphism of triangulated categories with duality in \cite[Definition 4.8]{handbook}. He specifies that $F$ commute with duality via an isomorphism $\eta : F\dual_\cat{A} \to \dual_\cat{B} F$ such that
$$
\xymatrix{
F \ar[r]^{F\Phi_\cat{A}} \ar[d]_{\Phi_\cat{B} F} & F \dual_\cat{A}^2 \ar[d]^{\eta \dual_\cat{A}}\\
\dual_\cat{B}^2 F \ar[r]_{\dual_\cat{B} \eta} & \dual_\cat{B} F \dual_\cat{A} 
}
$$
commutes. Given such an $\eta$ we can define a symmetric isomorphism 
$$
F \stackrel{F\Phi_\cat{A}}{\longrightarrow}F\dual_\cat{A}^2 \stackrel{\eta\dual_\cat{A}}{\longrightarrow}  \dual_\cat{B} F \dual_\cat{A} 
$$
and, conversely, given a symmetric isomorphism we can define such an $\eta$.
\end{remark}
\begin{proposition}
\label{functoriality prop}
A symmetric isomorphism of functors $\Psi:F \to \shift^r\dual_{\cat{A},\cat{B}} F$ of dimension $r$ induces a map of Witt groups
$$
W^*(F) : W^*(\cat{A}) \to W^{*+r}(\cat{B}).
$$
If $\Psi$ can be expressed as the boundary of a symmetric natural transformation of dimension $r+1$, \ie it fits into a diagram in $\funct{\cat{A}}{\cat{B}}$ of the form of that in Lemma \ref{boundary lemma}, then the induced map of Witt groups is zero.
\end{proposition}
\begin{proof}
We consider the case $*=r=0$ since all the others are similar. Choose a symmetric isomorphism $\alpha : \ob{A} \to \dual_\cat{A}\ob{A}$ representing a class in $W^0(\cat{A})$. We can construct a diagram
$$
\xymatrix{
F\ob{A}\ar[d]_{\Phi_\cat{B} F\ob{A}} \ar[rr]^{F\alpha} \ar[drrrr]^{\Psi \ob{A}}&& F\dual_\cat{A}\ob{A} \ar[rr]^{\Psi\dual_\cat{A}\ob{A}} && \dual_\cat{B} F\dual_\cat{A}^2\ob{A}  \ar[rr]^{\dual_\cat{B} F \Phi_\cat{A}\ob{A}} && \dual_\cat{B} F \ob{A} \ar@{=}[d]\\
\dual_\cat{B}^2F\ob{A} \ar[rr]_{\dual_\cat{B}^2F\Phi_\cat{A}\ob{A}} && \dual_\cat{B}^2F\dual_\cat{A}^2\ob{A} \ar[rr]_{\dual_\cat{B}\Psi \dual_\cat{A} \ob{A}}&& \dual_\cat{B} F\dual_\cat{A}\ob{A}\ar[rr]_{\dual_\cat{B} F \alpha}\ar[u]_{\dual_\cat{B} F \dual_\cat{A} \alpha} && \dual_\cat{B} F \ob{A}
}
$$ 
whose top row is an isomorphism $F\ob{A} \to \dual_\cat{B} F\ob{A}$ and whose bottom row is its dual. The lower triangle commutes because $\Psi$ is symmetric from $F$ to $\dual_{\cat{A},\cat{B}} F$. The upper triangle commutes because $\Psi$ is a natural transformation and the righthand square commutes because $\alpha$ is symmetric. Hence we have constructed a symmetric isomorphism $F\ob{A} \to \dual_\cat{B} F\ob{A}$ representing a class in $W^{0}(\cat{B})$. This is independent of the choice of representative $\alpha$.

If the symmetric natural transformation is a boundary we can explicitly construct a diagram expressing the isomorphism $F\ob{A} \to \dual_\cat{B} F\ob{A}$ as a boundary, and hence representing zero in the Witt group of $\cat{A}$.
\end{proof}

 \begin{example}
 \label{symmetry of phi}
 For any triangulated category $\cat{A}$ with $\delta$-duality the natural transformation $\Phi_{\cat{A}} : \id_{\cat{A}} \to \dual_{\cat{A}}^2$ is a symmetric isomorphism of functors because
 \begin{eqnarray*}
 \left(\dual_{\cat{A},\cat{A}}\Phi_{\cat{A}} \circ \Phi_{\cat{A},\cat{A}}(\id_{\cat{A}})\right)(\ob{A}) & = & \dual_{\cat{A}}\left( \Phi_{\cat{A}}(\dual_{\cat{A}}\ob{A}) \right) \circ \Phi_{\cat{A}}(\dual_{\cat{A}}^2\ob{A}) \circ \Phi_{\cat{A}}(\ob{A}) \\
 &=& \id_{\cat{A}}(\dual_{\cat{A}}^2\ob{A}) \circ \Phi_{\cat{A}}(\ob{A})\\
 &=&  \Phi_{\cat{A}}(\ob{A})
 \end{eqnarray*}
 using the definitions (\ref{functor duality}) and (\ref{functor phi}) and the coherence relation (\ref{coherence}). Conversely if $\Phi_\cat{A}$ is symmetric then it satisfies the coherence relation (\ref{coherence}). In other words the coherence relation precisely encodes the symmetry of the natural transformation $\Phi_\cat{A}$. Yet another way of phrasing this is that the identity functor is symmetrically self-dual via $\Phi_\cat{A}$.
 \end{example} 
 \begin{lemma}
 \label{comp of symm functors}
 If one is given symmetric natural transformations $F \to \dual_{\cat{A},\cat{B}} F$ and $G \to \dual_{\cat{B},\cat{C}} G $ then they can be composed (in two ways, which agree) to obtain a symmetric natural transformation
 $$
 G\circ F \to \dual_{\cat{B},\cat{C}} G \circ \dual_{\cat{A},\cat{B}} F  \to \dual_{\cat{A},\cat{C}} (G\circ F)
 $$
 where, since 
 \begin{eqnarray*}
 \dual_{\cat{B},\cat{C}} G \circ \dual_{\cat{A},\cat{B}} F & = & \dual_\cat{C} \circ G \circ \dual_\cat{B}^2 \circ F \circ \dual_\cat{A}\\
 \textrm{ and }\ \dual_{\cat{A},\cat{C}} (G \circ F) &= &\dual_\cat{C} \circ G \circ F \circ \dual_\cat{A},
 \end{eqnarray*}
 the final arrow arises from $\id \to \dual_\cat{B}^2$. A similar remark holds for symmetric natural transformations of other dimensions.
 \end{lemma}
 \begin{proof}
 Exercise!
 \end{proof}
 \begin{remark}
It follows from Example \ref{symmetry of phi} and Lemma \ref{comp of symm functors}  that there is a category whose objects are triangulated categories with duality and whose morphisms are symmetrically self-dual functors between them. 
 \end{remark}

 \subsection{Exact triples and long exact sequences}
\label{exact triples}

\begin{definition}
Recall that a full triangulated subcategory $\cat{A} \subset \cat{B}$ is \emph{thick} if, whenever $\ob{B}\oplus\ob{B}' \in \cat{A}$, then $\ob{B},\ob{B}' \in \cat{A}$.
\end{definition}
The importance of this concept is that if $\cat{A}$ is a thick subcategory the quotient category $\cat{B}/\cat{A}$, which has the same objects as $\cat{B}$ but in which all morphisms in $\cat{A}$ become invertible, inherits a triangulated structure (with respect to which the natural functor $\cat{B} \to \cat{B}/\cat{A}$ is triangulated). In this situation we say that
\begin{equation}
\label{exact triple}
\cat{A} \to \cat{B} \to \cat{B}/\cat{A}
\end{equation}
is an \emph{exact triple of triangulated categories}. If, in addition, $\cat{B}$ has a duality which preserves the subcategory $\cat{A}$ then $\cat{A}$ and $\cat{B}/\cat{A}$ inherit dualities and the inclusion and quotient functors are naturally symmetrically self-dual. We say (\ref{exact triple}) is an \emph{exact triple of triangulated categories with duality}.

\begin{theorem}[{Balmer \cite[Theorem 4.14]{handbook}}]
\label{balmer les}
Suppose $\cat{A} \to \cat{B} \to \cat{B}/\cat{A}$ is an exact triple of triangulated categories with duality in each of which $2$ is invertible. Then there is a long exact sequence of Witt groups$$
\ldots \to W^r(\cat{A})\to W^r(\cat{B})\to W^r(\cat{B} / \cat{A})\to W^{r+1}(\cat{A})\to \ldots
$$
\end{theorem}

\section{Witt groups of sheaves on a topological space}
\label{perfect witt groups}
Throughout this section all topological spaces will be assumed to be
Hausdorff, locally compact, locally connected and countable at infinity. We consider Witt groups of sheaves of $R$-modules on such spaces, where $R$ is a commutative regular Noetherian ring of finite Krull dimension in which $2$ is invertible.

\subsection{Perfect complexes of sheaves}

The bounded derived category $\der{X}$ of sheaves of $R$-modules over a space $X$ is a triangulated category. The left derived functor $\tensor$ of the tensor product of sheaves (obtained by taking flat resolutions) makes $\der{X}$ into a symmetric monoidal category. The unit  is the constant sheaf with stalk $R$, which we denote by $\struct{X}$. The right derived functor $\rshom{-}{-}$ of sheaf hom is an internal hom for this category. It is related to the morphisms in $\der{X}$ via taking the zero'th hypercohomology \ie
$$
\mor{\ob{E}}{\ob{F}} \cong \coh{0}{X}{\rshom {\ob{E}}{\ob{F}} }.
$$
These structures satisfy the properties described in \S\ref{internal structures}. Furthermore everything is enriched over the category of $R$-modules.

The bounded derived category comes equipped with a dualising object $\dualiser_X = p^!\struct{pt}$ where $p: X \to pt$ \cite[Chapter 3]{dimca}. However, although there is a natural transformation
\begin{equation}
\label{biduality}
\Phi_X : \id \to \rshom{\rshom{-}{\dualiser_X}}{\dualiser_X}
\end{equation}
it is not in general an isomorphism. Indeed, even over a point and when $R$ is a field we require finite dimensionality for a vector space to be isomorphic to its double dual. To fix this we pass to a subcategory. One choice of such a subcategory is the perfect derived category $\perf{X}$.
\begin{definition}[Verdier \cite{verdier}]
\label{perfect}
Let $x \in X$ and $U$ be a paracompact neighbourhood of $x$. Given $\ob{E} \in \der{X}$ let
$$
L(U; \ob{E}) \to \Gamma(U;\ob{I})
$$
be a projective resolution where $\ob{E} \to \ob{I}$ is an injective resolution of $\ob{E}$. We say $\ob{E}$ is \emph{perfect} if, given $U$, we can find $V \subset U$ such that the restriction factorises 
$$
\xymatrix{
L(U; \ob{E}) \ar[rr]^{\textrm{restriction}}  \ar@{-->}[dr]&& L(V; \ob{E})\\
&L_{U,V} \ar@{-->}[ur]
}
$$ 
via a bounded complex $L_{U,V}$ of projective modules of finite type.
\end{definition}
For the rings which we allow the following are equivalent \cite[Prop. 1.7]{verdier}:
\begin{enumerate}
\item $\ob{E}$ is perfect;
\item for any $x \in X$ and neighbourhood $U$ of $x$ there exists a smaller neighbourhood $V$ such that, for any $n \in \zz$, the restriction $H^n(U; \ob{E}) \to H^n(V; \ob{E})$ has a finitely generated image;
\item for any $x \in X$ and neighbourhood $U$ of $x$ there exists a smaller neighbourhood $V$ such that, for any $n \in \zz$, the extension $H_c^n(V; \ob{E}) \to H_c^n(U; \ob{E})$ has a finitely generated image.
\end{enumerate}
 
We recall some of the properties of perfect complexes.
\begin{proposition}[{Barthel \cite[\S 10]{sem-heid-stras}}]
The full subcategory $\perf{X}$ of perfect complexes is a triangulated subcategory of $\der{X}$. It contains the constant sheaf $\struct{X}$ and the dualising object $\dualiser_X$ and is preserved by $\rshom{-}{\dualiser_X}$. Furthermore the natural transformation $\Phi_X$ of (\ref{biduality}) becomes an isomorphism when restricted to $\perf{X}$.
\end{proposition}
\begin{proposition}[{Verdier \cite[Corollary 1.5]{verdier}}]
\label{perfect pushforward}
Suppose $f:X\to Y$ is a proper map. Then the derived pushforward $\rdf{f_*}:\der{X}\to\der{Y}$ takes perfect complexes to perfect complexes.
\end{proposition}

In order to obtain a theory which is functorial we will need to further restrict our attention to objects of $\perf{X}$ which have compact cohomological support in the following sense: 
\begin{definition} The \emph{cohomological support} of a bounded complex of sheaves is defined to be the union of the supports of its cohomology sheaves. This is clearly a quasi-isomorphism invariant. 
\end{definition}
It is an easy exercise to check that the full subcategory $\perfc{X}$ of objects with compact cohomological support is triangulated and preserved by duality. In addition, the natural transformation $\rdf{f_!} \to \rdf{f_*}$ becomes an isomorphism when restricted to objects with compact cohomological support and, by Proposition \ref{perfect pushforward}, we obtain a functor 
$$\rdf{f_!} = \rdf{f_*} : \perfc{X} \to \perfc{Y}.$$
\begin{definition}
The \emph{perfect Witt group} $W_r^p(X)$ of  $X$ is defined to be the Witt group $W^{-r}(\perfc{X})$. Note the switch to homological indexing reflecting the fact that the perfect Witt groups will turn out to be covariant functors.
\end{definition}

%Later, in \S\ref{excision}, we will require a subcategory with even better properties than the perfect %derived category. To find such a subcategory we will restrict $X$ to be a stratifiable space and %consider the constructible derived category $\constr{X}$.

\subsection{Functoriality}

What are the functorial properties of the perfect Witt groups $W^p_*(X)$? To answer this question we study the relationship between Verdier dual pairings and the pushforward with proper support $\rdf{f_!}$ where $f: X \to Y$ is a continuous map. First we give a criterion for a map $\ob{E} \to \ob{F}$ to be an isomorphism in $\der{X}$.
\begin{lemma}
\label{IM criterion}
A map $\ob{E} \to \ob{F}$ is an isomorphism in $\der{X}$ if, and only if, for all open $U$ in $X$ the induced map $\rdf{p_*}\imath^*\ob{E} \to \rdf{p_*}\imath^*\ob{F} $ in $\der{\mathrm{pt}}$ is an isomorphism, where
$$
\mathrm{pt} \stackrel{p}{\longleftarrow} U \stackrel{\imath}{\hookrightarrow} X.
$$
\end{lemma}
\begin{proof}
This follows from the fact that a map in $\der{X}$ is an isomorphism if, and only if, it induces isomorphisms on all cohomology sheaves.
\end{proof}

\begin{lemma}
\label{vdp criterion}
Given a map $\alpha : \ob{E} \tensor \ob{F} \to \dualiser_X$ there is an induced map $$
\rdf{p_*}\imath^* \ob{E} \tensor \rdf{p_!}\imath^* \ob{F} \to \rdf{p_!}\imath^*\left( \ob{E} \tensor \ob{F} \right) \to \rdf{p_!}\imath^*\dualiser_X \to \dualiser_{\mathrm{pt}}
$$
in $\perf{\mathrm{pt}}$ for each open $U$ in $X$, where $\mathrm{pt} \stackrel{p}{\longleftarrow} U \stackrel{\imath}{\hookrightarrow} X$ and the final map arises from the identification $\imath^*\dualiser_X = \dualiser_U = p^!\dualiser_ {\mathrm{pt}}$ and the unit $\rdf{p_!}p^! \to \id$ of the adjunction. The map $\alpha$ is a Verdier dual pairing, \ie induces an isomorphism $\ob{E} \to \dual\ob{F}$, if, and only if, the induced map is a Verdier dual pairing for every open $U$.
\end{lemma}
\begin{proof}
This follows directly from the criterion in Lemma \ref{IM criterion}.
\end{proof}
\begin{lemma}
\label{vdp pushforward}
Given  a Verdier dual pairing $\alpha : \ob{E} \tensor \ob{F} \to \dualiser_X$ and a map $f:X \to Y$ there is an induced Verdier dual pairing
$$
\rdf{f_*} \ob{E} \tensor \rdf{f_!} \ob{F} \to \rdf{f_!}\left( \ob{E} \tensor \ob{F} \right) \to \rdf{f_!}\dualiser_X \to\dualiser_Y.
$$
\end{lemma}
\begin{proof}
This follows directly from Lemma \ref{vdp criterion} since restriction to an open commutes with both $\rdf{f_*}$ and $\rdf{f_!}$.
\end{proof}
\begin{proposition}
\label{pushforward}
For any $f:X \to Y$ there is a symmetric natural transformation $\chi_f:\rdf{f_!} \to \dual \rdf{f_!}\dual$ of functors $\perf{X} \to \perf{Y}$ which is an isomorphism when restricted to objects in $\perfc{X}$. This induces maps
$$
f_* : W^p_*(X) \to W^p_*(Y)
$$
satisfying $(f\circ g)_* = f_* \circ g_*$ and $\id_* = \id$.
\end{proposition}
\begin{proof}
Given $\ob{E} \in \perf{X}$ there is a natural Verdier dual pairing $\ob{E} \tensor \dual\ob{E} \to \dualiser_X$ corresponding to the isomorphism $\ob{E} \to \dual^2\ob{E}$. By Lemma \ref{vdp pushforward} this yields a Verdier dual pairing
$$
\rdf{f_*} \ob{E} \tensor \rdf{f_!} \dual\ob{E} \to \dualiser_Y
$$
corresponding to a natural isomorphism $\rdf{f_*} \ob{E} \to \dual\rdf{f_!} \dual \ob{E}$. We define the required natural transformation $\chi_f$ on $\ob{E}$ by precomposing with $\rdf{f_!}\ob{E} \to \rdf{f_*}\ob{E}$, which is an isomorphism if $\ob{E}$ has compact cohomological support. Symmetry follows from the symmetry of the natural transformation $\id \to \dual^2$ (see Example \ref{symmetry of phi}). 

The induced map $f_*: W^p_*(X) \to W^p_*(Y)$ is defined as in \S \ref{functors and duality}. Its properties follow easily.
\end{proof}
\begin{remark}
If we consider the Witt groups of $\perf{X}$ rather than $\perfc{X}$ we obtain a theory `with closed supports' which is functorial under \emph{proper} maps.
\end{remark}

\subsection{Homotopy invariance}
\label{htpy}
\begin{proposition}
Suppose $h:X \times I \to Y$ is a homotopy between $f:X \to Y$ and $g: X \to Y$. Then the induced maps $f_*$ and $g_*$ from $W_*^p(X)$ to $W_*^p(Y)$ agree.
\end{proposition}
\begin{proof}
We have maps
$$
\xymatrix{
X \times (0,1)  \ar@{^{(}->}[r]^{\imath} \ar[dr]_q & X \times [0,1] \ar[d]^p& X  \ar@{_{(}->}[l]_{\jmath} \times \{0,1\} \\
& X
}
$$
Since $q$ is a projection with smooth fibre of dimension $1$ we have $q^*\dualiser_X \cong q^!\dualiser_X \cong \shift \dualiser_{X \times(0,1)}$. Starting with the standard Verdier dual pairing $\ob{E} \tensor \dual\ob{E} \to \dualiser_X$ for $\ob{E} \in \perfc{X}$ we obtain a map
$$
q^*\ob{E} \tensor q^*\dual \ob{E} \cong q^*\left( \ob{E} \tensor \ob{E} \right) \to q^*\dualiser_X \cong \shift \dualiser_{X \times(0,1)}.
$$
This defines a symmetric natural transformation 
$q^* \to \shift^{-1}\dual q^*\dual$
(of dimension $-1$) of functors from $\perfc{X}$ to $\der{X\times(0,1)}$. It is easy to check that $q^*$ preserves perfect complexes so that the image is in the full triangulated subcategory $\perf{X\times(0,1)}$. It follows from Proposition \ref{pushforward} that we also have a symmetric natural transformation $\chi_{\imath}: \rdf{\imath_!} \to \dual \rdf{\imath_!} \dual$ of functors from $\perf{X\times(0,1)}$ to $\perf{X\times[0,1]}$. By Lemma \ref{comp of symm functors} we can compose the two to obtain a symmetric natural transformation 
$$
\rdf{\imath_!}q^* \to \shift^{-1}\dual\rdf{\imath_!}q^*\dual
$$
of dimension $-1$. Note that this will not be an isomorphism, but that we can explicitly identify the coboundary to be the natural isomorphism $\rdf{\jmath_*} \jmath^* q^* \to \dual \left( \rdf{\jmath_*} \jmath^* q^*\right)\dual$ given by the matrix $$\left(\begin{matrix} \Phi_X & 0 \\ 0 & -\Phi_X \end{matrix}\right)$$ with respect to the decomposition corresponding to the two components of $X \times \{0,1\}$.  Composing with the symmetric natural isomorphism $\chi_h$ we exhibit $\chi_f \oplus (-\chi_g)$ as a coboundary. Hence, by the last part of Proposition \ref{functoriality prop}, we have $f_*-g_*=0$.
\end{proof}

\subsection{Relative Witt groups and long exact sequences}
\label{relative groups}
Suppose $\jmath : A \hookrightarrow X$ is a closed inclusion. Then $\rdf{\jmath_!} : \perfc{A} \to \perfc{X}$ is the inclusion of a thick subcategory. It follows from \S\ref{exact triples} that we have an exact triple of triangulated categories with duality
$$
\perfc{A} \to \perfc{X} \to  \perfc{X}/ \perfc{A}.
$$
\begin{definition}
We define the {\em relative Witt group} $W_r^p(X,A)$ to be the $(-r)^{th}$ Witt group of $ \perfc{X}/ \perfc{A}$. 
\end{definition}
An immediate consequence of Balmer's theorem (\ref{balmer les}) is that we obtain a long exact sequence
\begin{equation}
\label{les}
\ldots \to W^p_r(A) \to W^p_r(X) \to W^p_r(X,A) \to W^p_{r-1}(A) \to \ldots
\end{equation}
\begin{remark}
The functoriality and homotopy invariance of the perfect Witt groups  can be extended to the relative groups.
\end{remark}

\subsection{Constructibility and excision} 
\label{excision}
The axiom of a generalised homology theory with which we have not yet dealt is excision. We would like to know that if $\overline U \subset A^\circ \subset X$ (with $U$ open and $A$ closed) then the map
$$
W^p_*(X - U, A - U) \to W^p_*(X , A)
$$
induced by inclusion of pairs is an isomorphism. The simplest proof would be to show that the functor
$$
\perfc{X - U} / \perfc{A - U} \to \perfc{X} / \perfc{A}
$$
induced from extension by zero is an equivalence. Unfortunately, this does not seem easy (and indeed may not be true without further assumptions). This is the same difficulty confronted in Hutt's paper \cite{hutt} and we resolve it in a similar way, by restricting both the types of space and of sheaves that we consider.

Henceforth all our spaces and maps will be (topologically) stratified in the sense of \cite[\S 1]{gm2}. A space $X$ with a given stratification $S$ will be denoted by $X_S$. Let $\constr{X_S}$ denote the full triangulated subcategory of $\der{X}$ consisting of complexes of sheaves whose cohomology is constructible with respect to $S$ (in the sense of \cite[\S 1]{gm2}). Recall from \cite[\S 1]{gm2} that for any stratified space $X_S$
\begin{enumerate}
\item  $\constr{X_S}$ is triangulated;
\item $\constr{X_S}\subset \perf{X}$;
\item the dualising object $\dualiser_X \in \constr{X_S}$;
\item duality preserves  $\constr{X_S}$.
\end{enumerate}
It follows that both $\constr{X_S}$ and the full subcategory  $\constrc{X_S}$ of objects with compact cohomological support are triangulated categories with duality. 

\begin{definition}
For a stratified space $X_S$ we define the {\em constructible Witt groups} $W^c_*(X_S) = W_*(\constrc{X_S})$.
\end{definition}

If $f:X_S \to Y_T$ is a stratified map then it induces a functor $\rdf{f_*} :  \constr{X_S}\to  \constr{Y_T}$ with a left adjoint $f^* : \constr{Y_T}\to \constr{X_S}$ and, dually, $\rdf{f_!}$ with a  right adjoint $f^!$. The natural map $\rdf{f_!} \to \rdf{f_*}$ becomes an isomorphism on $\constrc{X_S}$. For the special case of a complementary pair of open and closed inclusions
$$
U_S \stackrel{\imath}{\hookrightarrow} X_S \stackrel{\jmath}{\hookleftarrow} A_S
$$
we obtain `glueing data':
\[
\xymatrix{
\constr{U_S} \ar@<3pt>[rr]^{\rdf{\imath_*}\ ,\ \rdf{\imath_!}} &&\constr{X_S} \ar@<3pt>[rr]^{ \jmath^*\  ,\ \jmath^!} \ar@<3pt>[ll]^{\imath^*=\imath^!} &&\constr{A_S} \ar@<5pt>[ll]^{\rdf{\jmath_!} = \rdf{\jmath_*}}
}
\]
obeying the usual relations, see, for example, \cite[Chapter 5 \S 3.9.1]{gema}.

By making small modifications to our previous arguments we can show that the constructible Witt  groups are functorial under stratified maps, and are stratified homotopy invariants. (Stratified maps $f,g : X_S \to Y_T$ are \emph{stratified homotopic} if there is a homotopy $h$ from $f$ to $g$ which is a stratified map
$$
h: X \times [0,1] \to Y_T
$$
where the stratification of $X \times [0,1]$ is the product of the given stratification of $X$ and the obvious stratification of $[0,1]$ by the interior and endpoints.) Closed stratified inclusions induce long exact sequences involving the relative groups $W^c_*(X,A) = W_*(\constrc{X}/\constrc{A})$. 

\begin{lemma}
If $\overline U_S \subset A_S^\circ \subset X_S$ are stratified inclusions (with $U$ open and $A$ closed) then we have an equivalence
$$
\xymatrix{
\frac{\constrc{X_S - U_S} }{\constrc{A_S - U_S}}\ar@<3pt>[r]^{\rdf{\jmath_*}} &\frac{\constrc{X_S }}{\constrc{A_S }} \ar@<3pt>[l]^{\jmath^*} \\
}
$$
where $\jmath : X - U \hookrightarrow X$.
\end{lemma}
\begin{proof}
We have $\jmath^*\rdf{\jmath_*} \cong \id$ and there is a triangle 
$$
\rdf{\imath_! \imath^!} \to \id \to \rdf{\jmath_*}\jmath^* \to \shift \rdf{\imath_! \imath^!} 
$$
in $\constrc{X}$ where $\imath : U \hookrightarrow X$. Since $\overline U \subset A^\circ$ we see that $\rdf{\imath_! \imath^!}  \cong 0$ as an endofunctor of $\constrc{X_S } / \constrc{A_S }$ and so $\rdf{\jmath_*}\jmath^* \cong \id$ too.
\end{proof}
We immediately obtain
\begin{corollary}
\label{constructible excision}
The constructible Witt groups satisfy excision for stratified maps.
\end{corollary}

\subsection{Constructible Witt groups of polyhedra}
\label{constructible witt groups}
Let $K$ be  a simplicial complex. Its realisation is naturally stratified (see \S\ref{comb and constr}), and we denote this stratified space by $K_S$. (The underlying topological space is just the realisation $|K|$ but we want to remember the stratification.) We write $\constr{K_S}$ for the constructible derived category of sheaves on $K_S$. The realisation of a simplicial map $K\to L$ is a stratified map with respect to these natural stratifications. Furthermore
\begin{lemma}
If simplicial maps $f,g : K \to L$ are contiguous then their realisations $|f|$ and $|g|$ are stratified homotopic.
\end{lemma}
\begin{proof}
We can reduce to the case when $K$ consists of a single simplex. Interpolating linearly between $f$ and $g$ then gives the desired stratified homotopy.
\end{proof}
Since the constructible Witt groups are stratified homotopy invariant functors (see \S \ref{excision}) it follows from this lemma that we obtain combinatorial homotopy invariant functors $K \mapsto W^c_i(K_S)$ from the category of simplicial complexes and simplicial maps.
\begin{theorem}
\label{simplicial ght}
For a commutative regular Noetherian ring $R$ of finite Krull dimension in which $2$ is invertible the constructible Witt groups form a combinatorial generalised homology theory on simplicial complexes and maps.
\end{theorem}
 \begin{proof}
We have already seen that the constructible Witt groups define combinatorial homotopy invariant functors. It remains to check that they satisfy excision and that there is a relative long exact sequence associated to any pair. Excision follows from Corollary \ref{constructible excision} since the realisation of a simplicial map is stratified. The long exact sequences for pairs arise as in \S \ref{relative groups}. 
 \end{proof}
 
Note that whenever $K'$ is a refinement of a simplicial complex $K$ then there is an induced inclusion $\constr{K_S} \hookrightarrow \constr{K'_S}$ of triangulated categories  with duality. 
\begin{definition}
Suppose $X$ is a compact polyhedron \ie a compact Hausdorff topological space with a chosen family of compatible triangulations, any two of which have a common refinement. We define the {\em constructible derived category}
$$
\constrc{X} = \mathrm{colim}_{|K|=X} \constr{K_S}.
$$
This is a triangulated category and inherits a duality. (It is unfortunate, but seemingly unavoidable, that the term `triangulated' is used here in both the geometric and categorical senses.) We denote the Witt groups of this category  by $W^c_*(X)$ and call them the constructible Witt groups of $X$. They are independent of any particular triangulation. In fact
\end{definition}
\begin{theorem}
\label{ght}
For a commutative regular Noetherian ring $R$ of finite Krull dimension in which $2$ is invertible the constructible Witt groups $W^c_i(-)$ form a generalised homology theory on compact polyhedra and continuous maps. In particular they are homotopy invariant functors. 
\end{theorem}
 \begin{proof}
First of all we need to consider functoriality. Suppose $f: X \to Y$ is a continuous map of compact polyhedra. Given any triangulations $K$ and $L$ of $X$ and $Y$ respectively the simplicial approximation theorem says that we can find a refinement $K'$ and a simplicial approximation $K' \to L$ of $f$. Furthermore any two such approximations are combinatorially homotopic. Since the constructible Witt groups of simplicial complexes are combinatorial homotopy invariant functors it follows that the constructible Witt groups of polyhedra are homotopy invariant functors.

Excision and the existence of long exact sequences for pairs now follow from the combinatorial analogues.
 \end{proof}
 \begin{remark}
 \label{ght remark}
Note that it follows that if $K$ is a triangulation of the compact polyhedron $X$ then the inclusion $\constrc{K_S} \hookrightarrow \constrc{X}$ induces isomorphisms $W_i^c(K_S) \cong W^c_i(X)$ for all $i$.
\end{remark}

\begin{remark}
Since constructible sheaves are perfect $\constrc{X}$ includes in $\perfc{X}$ and thence there are maps $W^c_*(X) \to W^p_*(X)$. The extent to which these fail to be isomorphisms measures the failure of the perfect Witt groups to be a generalised homology theory. Indeed, since the direct sum of two sheaves is locally constant if, and only if, both summands are locally constant, we see that $\constr{K}$ is a thick subcategory of $\perfc{X}$. It follows from Theorem \ref{balmer les} that there is a long exact sequence
$$
\ldots \to W^c_i(X) \to W^p_i(X) \to W^{-i}(\perfc{X} / \constrc{X}) \to W^c_{i-1}(X) \to \ldots
$$
relating the constructible and perfect Witt groups of a polyhedron. Unfortunately I do not know of an example in which the relative term is non-zero.
\end{remark}

\subsection{Products}
\label{products}

In this section we assume that $R$ is a field. For any stratified space, and thence for a polyhedron $X$, the derived tensor product $\tensor$ preserves constructible complexes. Thus given polyhedra $X$ and $Y$ we can define an external tensor product
\begin{eqnarray*}
\boxtimes : \constrc{X} \times \constrc{Y} &\to& \constrc{X \times Y}\\
(\ob{E},\ob{F}) &\mapsto& p^*\ob{E} \tensor q^*\ob{F}
\end{eqnarray*}
where $p$ and $q$ are the projections onto $X$ and $Y$ respectively.

The product category $\constrc{X} \times \constrc{Y} $ inherits a product duality given by $(\ob{E},\ob{F}) \mapsto (\dual_X\ob{E},\dual_Y\ob{F})$. The dualising object $\dualiser_{X \times Y} \cong \dualiser_X\boxtimes \dualiser_Y$ (this is essentially the Kunneth theorem --- see \cite[\S 1.12]{gm2} or \cite[p181]{bo}). Hence, using \cite[Corollary 2.0.4]{MR2031639}, we have
\begin{eqnarray*}
\dual (\ob{E} \boxtimes \ob{F}) & \cong &  \rshom{\ob{E} \boxtimes\ob{F}}{\dualiser_{X \times Y}} \\
& \cong &  \rshom{\ob{E} \boxtimes\ob{F}}{\dualiser_X\boxtimes \dualiser_Y} \\
& \cong &  \rshom{\ob{E} }{\dualiser_X} \boxtimes  \rshom{\ob{F} }{\dualiser_Y}\\
& \cong &  \dual\ob{E} \boxtimes \dual\ob{F}.
\end{eqnarray*}
In other words $\boxtimes$ intertwines the product duality with the standard Poincar\'e--Verdier duality on $\constrc{X \times Y}$. It follows that the external tensor product is a map of categories with duality and so induces maps
$$
W_i^c(X) \times W_j^c(Y) \to W_{i+j}^c(X \times Y).
$$
These are easily seen to factorise through the tensor product so that we have a graded product
$$
W_*^c(X) \otimes W_*^c(Y) \to W_{*}^c(X \times Y).
$$
In particular $W^c_*(X)$ is always a $W^c_*(\textrm{pt})$-module. This product can be extended to relative groups.

\subsection{Generalisations and related theories}
\label{other theories}

We briefly touch upon some of the possible generalisations and other theories which can be constructed using similar techniques. The first remark is that everything can be done equivariantly. Suppose that a group $G$ acts (piecewise-linearly) upon a polyhedron $X$. Then there is an equivariant constructible derived category $\eqconstr{G}{X}$ of sheaves on $X$ (see \cite{bl}), equipped with an equivariant Poincar\'e--Verdier duality. The $G$-equivariant constructible Witt groups of $X$ are the Witt groups of this category. Rewriting our previous arguments, using the technology of functors between equivariantly constructible derived categories of sheaves developed in \cite{bl}, we see that these are functorial for equivariant maps, are equivariant homotopy invariants etc. (The hard work is in the construction of the equivariant derived category, once we have that the definition and properties of the Witt groups are routine.)

In a different direction, we can define new theories by restricting the allowed stratifications. For a complex algebraic variety $V$ the Witt groups of the derived category $\constrac{V}$ of sheaves whose cohomology is constructible with respect to a stratification of $V$ by complex algebraic varieties are considered in \cite{yokura,chern-classes}. The Riemann--Hilbert correspondence provides us with an alternative description for these groups as the Witt groups of the derived category of regular holonomic $\mathcal{D}$-modules on $V$. 

One can also study the constructible Witt groups $W_*^c(X_S)$ of a
stratified space (with fixed stratification) in their own right,
rather than using them as a tool to obtain a stratification-invariant
theory as we have done. This is the approach taken in \cite{cs}, where
they obtain a powerful `decomposition theorem up to cobordism'
\cite[Theorem 4.7]{cs} by identifying a set of generators of their
cobordism group which  correspond to irreducible self-dual perverse sheaves supported on the strata. 

\section{L-theory and constructible Witt groups}
\label{L-theory}

We relate the constructible Witt groups of a simplicial complex to the free symmetric L-groups, showing that,  under certain conditions, they are isomorphic. 

In this section $R$ will be a commutative ring and $K$ a finite simplicial complex.

\subsection{$\rx$-modules and L-theory}

In \cite{ranweiss} Ranicki and Weiss define an $\rx$-module to be  a finitely generated free $R$-module $A$ with a direct sum decomposition
$$
A =  \bigoplus_{\sigma \in K} A(\sigma)
$$
into free $R$-modules. A map of $\rx$-modules is an $R$-module morphism $\alpha : A \to B$ such that 
$$
\alpha(A(\sigma)) \subset \bigoplus_{\tau \geq \sigma} B(\tau).
$$
(Here we regard $\sigma \leq \tau$ if $\sigma$ is a face of $\tau$.)
More concisely, if we regard $K$ as a category with objects the simplices and morphisms the face inclusions, then an $\rx$-module is a functor from $K^{op}$ to $R$-modules and a map of $\rx$-modules is a natural transformation. The category $\rxmod$ of $\rx$-modules is a full subcategory of the functor category $\funct{K^{op}}{\rmod}$.  

\begin{example}
Define a chain complex $\cochains{}{X}$ of $\rx$-modules  by 
$$
\cochains{}{X}_{i}(\sigma) = \left\{ 
\begin{array}{ll}
R & i=-\dim \sigma\\
0 & i\neq -\dim \sigma.
\end{array}
\right.
$$
Thus $\cochains{}{X}(U)$ is the group of simplicial cochains on $U$ re-graded as a chain complex. The differentials are the coboundary maps. 
\end{example}

The category $\rxmod$ is additive and has a a natural chain duality in the sense of \cite[Definition 1.1]{ranicki} \ie a functor $T$ which takes an $\rx$-module to a bounded complex of $\rx$-modules: 
$$
T(A) = \hom{R}{\hom{(R,K)}{\cochains{}{X}}{A}}{R}
$$
(with the $\rx$-module structure given by
$$
T(A)_i(\sigma) = \left\{ 
\begin{array}{ll}
\bigoplus_{\tau \geq \sigma} \hom{R}{A(\tau)}{R} & i = -\dim \sigma \\
0 & i \neq -\dim \sigma)
\end{array}
\right.
$$ 
and a natural transformation $e:T^2 \to 1$ satisfying 
\begin{enumerate}
\item $e(TA)\circ (T(e(A)) = 1$ and
\item $e(A)$ is a chain equivalence. 
\end{enumerate}
\begin{example}
\label{dual example}
For $\sigma \in K$ define an $\rx$-module $C_\sigma$ by
$$
C_\sigma(\tau) = \left\{
\begin{array}{ll}
R & \tau = \sigma \\
0 & \tau \neq \sigma
\end{array}
\right.
$$
Then $TC_\sigma$ is the chain complex of $\rx$-modules with
$$
(TC_\sigma)_i(\tau) = \left\{
\begin{array}{ll}
R & \tau \leq \sigma\ \textrm{and}\  i=-\dim \tau \\
0 & \textrm{otherwise.}
\end{array}
\right.
$$
and differentials given by the coboundary maps.
\end{example}
We will write $\comrxmod$ for the category of bounded chain complexes of $\rx$-modules and chain maps. Since we are working with chain complexes, and not cochain complexes as elsewhere in this paper,  $\Sigma$ will denote the right shift functor \ie for $A \in \comrxmod$ we have $(\Sigma A)_i = A_{i-1}$. The chain duality $T$ extends to a duality on $\comrxmod$ with $\Sigma \circ T = T \circ \Sigma^{-1}$ (see \cite[p.26]{ranicki}). 
\begin{definition}
\label{spc}
An {\em $n$-dimensional symmetric Poincar\'e complex} is an element $A \in \comrxmod$ together  with a collection of maps
$$
\{ \phi_s \in \mor{TA}{\Sigma^{s-n}A} : s \geq 0\} 
$$
such that $\phi_0$ is an $\rx$-module chain equivalence (\ie $\phi_0(\sigma)$ is a quasi-iso\-morphism of chain complexes of $R$-modules for each $\sigma$) and, for $s \geq 0$, 
$$
\phi_s + (-1)^{s+1} e(\Sigma^{s-n}A) \circ \Sigma^{s-n}T ( \phi_s ) = \partial_{\Sigma^{s-n}A} \circ \phi_{s+1} +   \phi_{s+1}\circ \partial_{TA}. 
$$
The idea is that $\phi_0$ is an $\rx$-module chain equivalence which is symmetric up to a homotopy $\phi_1$, which in turn is symmetric up to a homotopy $\phi_2$, and so on. 
\end{definition}
\begin{definition} The {\em free symmetric L-group} $H_n(K; \mathbb{L}^\cdot(R))$ is an Abelian group generated by the $n$-dimensional symmetric Poincar\'e complexes modulo the algebraic cobordism relation defined in \cite[Definition 1.7]{ranicki}. It is the homology group of the symmetric L-theory spectrum $\mathbb{L}^\cdot(R)$ in \cite[Proposition 14.5]{ranicki}. 
\end{definition}

A simplicial map $f:K \to L$ induces a pushforward $f_*$ from $\rxmod$ to $\textrm{(R,L)-Mod}$ given by 
$$
(f_*A)(\sigma) = \bigoplus_{f(\tau) = \sigma} A(\tau)
$$ 
or, equivalently considering $A \in \funct{K^{op}}{\rmod}$, by composing with $f^{op} : K^{op} \to L^{op}$. This yields a map $f_*: H_*(K; \mathbb{L}^\cdot(R)) \to H_*(L; \mathbb{L}^\cdot(R))$, making the free symmetric L-groups functorial. In fact they form a generalised homology theory --- see \cite[\S 12 -- 14]{ranicki}. 

\subsection{Combinatorial and constructible sheaves}
\label{comb and constr}
In order to relate the free symmetric L-groups to the constructible Witt groups we need to establish a relationship between $\rx$-modules and constructible sheaves. This is done by relating both to combinatorial sheaves on $K$, which we now describe.

 We can topologise $K$ by defining a subset of simplices to be open if, and only if, it is upwardly closed \ie a set $U$ is open if, and only if, $\tau \in U$ whenever there is some $\sigma \leq \tau$ with $\sigma \in U$.  There is a unique smallest open set containing any simplex $\sigma$, namely the star $\st{\sigma} = \bigcup_{\tau \geq \sigma} \tau$. The stars form a base for this topology. 
 
 We will use the term  combinatorial sheaf on $K$ to describe  a sheaf of finitely generated $R$-modules in this topology and denote the category of combinatorial sheaves by $\sh{K}$.
\begin{lemma}
Combinatorial sheaves  are precisely functors from $K$ to the category of finitely generated $R$-modules. Maps of combinatorial sheaves are natural transformations.
\end{lemma}
\begin{proof}
A combinatorial sheaf $\mathcal{E}$ is equivalent to the following data: a finitely generated  $R$-module $\mathcal{E}(\st{\sigma})$ for each basic open set $\st{\sigma}$ and a set of compatible maps $\mathcal{E}(\st{\sigma})\to \mathcal{E}(\st{\tau})$, one for each inclusion $\st{\tau} \subset \st{\sigma}$. Since $\st{\tau} \subset \st{\sigma}$ if, and only if, $\tau \geq \sigma$ it is clear that the assignment
$$
\sigma \mapsto \mathcal{E}(\st{\sigma})
$$
defines a functor. Conversely, a functor defines a combinatorial sheaf. The final statement follows easily.
\end{proof}

The geometric realisation $|K|$ of $K$ has a natural stratification with strata the $S_\sigma$ where $S_\sigma = |\sigma| - |\partial \sigma|$. As in \S\ref{constructible witt groups}, when we want to emphasize that $|K|$ is a stratified space we denote it by $K_S$. We say a sheaf $\mathcal{E}$ of finitely generated $R$-modules on $K_S$ is constructible if the restriction $\mathcal{E}|_{S_\sigma}$ is constant for each $\sigma \in K$ and denote the subcategory of constructible sheaves of finitely generated $R$-modules by $\shc{K_S}$.

There is a continuous map $s:K_S \to K$ given by $x \mapsto \sigma$ where $x \in S_\sigma$. \begin{lemma}
\label{constr sheaves lemma}
The category of combinatorial sheaves on $K$ is equivalent, indeed isomorphic, to the category of constructible sheaves on $K_S$ via:
$$
\xymatrix{
\shc{K_S} \ar@<3pt>[r]^{s_*} &\sh{K} \ar@<3pt>[l]^{s^*} 
}
$$
\end{lemma}
\begin{proof}
For $\mathcal{E} \in \shc{K_S}$ we have $(s_*\mathcal{E})_\sigma = s_*\mathcal{E}(\st{\sigma}) = \mathcal{E}(|\st{\sigma}|)$ and for  $\mathcal{F} \in \sh{K}$ we have $(s^*\mathcal{F})_x = \mathcal{F}_{sx} = \mathcal{F}_\sigma = \mathcal{F}(\st{\sigma})$ where $x \in S_\sigma$. In particular, note that $s^*\mathcal{F}$ is a constructible sheaf. Hence we have, on the one hand, $(s_*s^*\mathcal{F})_\sigma = s^*\mathcal{F}(|\st{\sigma}|) = \mathcal{F}(\st{\sigma}) = \mathcal{F}_\sigma$ and, on the other, $(s^*s_*\mathcal{E})_x = (s_*\mathcal{E})_\sigma = \mathcal{E}(|\st{\sigma}|)$ where $x \in S_\sigma$. Now there is clearly a map
$$
\mathcal{E}(|\st{\sigma}|) \to  \mathcal{E}_x
$$
and, since $|\st{\sigma}|$ is a union of contractible strata on each of which $\mathcal{E}$ is constant, it is an isomorphism. Hence $s_*s^*$ and $s^*s_*$ are both the identity.
\end{proof}
The pullback $s^*$ extends to a triangulated functor 
\begin{equation}
\label{pullback}
s^* : \der{K} \to \constr{K_S}
\end{equation}
where $\der{K}$ is the bounded derived category of $\sh{K}$. Note that we do not need, or claim, that this is an equivalence because $\constr{K_S}$ consists of cohomologically constructible complexes and is not necessarily equivalent to $\der{\shc{K_S}}$.

%\begin{corollary}
%The bounded derived category $\der{K}$ of $\sh{K}$ is equivalent to the bounded constructible %derived category $\constr{K_S}$. 
%\end{corollary}
%\begin{proof}
%This is not quite so obvious as it might at first seem since 
%
%INDEED I DOUBT IT'S TRUE!
%
%objects of the constructible derived category 
%need not be bounded complexes of constructible sheaves, they are bounded complexes whose %cohomology sheaves are constructible i.e. in $\shc{K_S}$. However there is clearly a functor
%$$
%s^*:\der{K} \to \constr{K_S}.
%$$
%Both categories have bounded t-structures, $s^*$ is t-exact and induces an isomorphism between the %hearts by Lemma \ref{constr sheaves lemma}. It follows that $s^*$ induces an equivalence --- see \eg %\cite[Chapter 5, 3.7.3]{gema}.
%\end{proof}

\subsection{From L-groups to Witt groups}
We now explain how to construct a natural transformation from the free symmetric L-groups to the constructible Witt groups.

\begin{lemma}
\label{fff lemma}
There is a fully faithful  functor from $\rxmod^{\, op}$ to $\sh{K}$ which takes an $\rx$-module $A$ to the combinatorial sheaf $\mathcal{A}$ with $$\mathcal{A}(U) =  \bigoplus_{\sigma \in U} \mor{A(\sigma)}{R}$$ (with restriction maps given by the obvious projections). Furthermore this functor is natural in the sense that it takes the pushforward $f_*A$ of an $\rx$-module under a simplicial map $f$ to the pushforward $f_*\mathcal{A}$ of the corresponding combinatorial sheaf. 
\end{lemma}
\begin{proof}
This follows immediately from the above description of the category of combinatorial sheaves as functors from $K$ to $\rmod$ and the fact that $\rx$-modules form a full subcategory of $\funct{K^{op}}{\rmod}$.
\end{proof}

Clearly we can extend this to a functor taking chain complexes of $\rx$-modules to cochain complexes of combinatorial sheaves. Composing this with the functor (\ref{pullback}) from $\der{K}$ to $\constr{K_S}$ we obtain a functor
$$
F: \comrxmod^{op} \to \constr{K_S}.
$$

Suppose that $R$ is regular, Noetherian and of finite Krull dimension so that there is a Poincar\'e--Verdier duality functor $D: \constr{K_S}^{op} \to \constr{K_S}$. Then we have
\begin{lemma}
\label{duality lemma}
The functor $F$  commutes with duality i.e. $F \circ T = D \circ F$.
\end{lemma}
\begin{proof}
The statement of the lemma needs clarification since the Poincar\'e--Verdier dual is only defined up to isomorphism in $\constr{K_S}$. Thus, in order to make sense of the statement, we need to specify representative complexes of sheaves for each $D\circ F(A)$. Note that it is sufficient to do this for each  $\rx$-module of
the form $C_\sigma$ with
$$
C_\sigma(\tau) = \left\{
\begin{array}{ll}
R & \tau = \sigma \\
0 & \tau \neq \sigma.
\end{array}
\right.
$$
and that it suffices to prove that $F \circ T(C_\sigma) =
D \circ F(\sigma)$ for each $\sigma \in K$.

We have $F(C_\sigma) = \jmath_*\jmath^* \mathcal{O}_{K_S}$ where $\jmath : |\bar{\sigma}| \hookrightarrow K_S$ is the inclusion and $\mathcal{O}_{K_S}$ the constant sheaf with stalk $R$ on $K_S$. Hence $D\circ F(C_\sigma)$ is isomorphic to the pushforward $\jmath^* \mathcal{D}_{|\overline{\sigma}|}$ of the dualising complex on $|\overline{\sigma}|$. We can choose to represent it by the complex $s^*\mathcal{C}$ of sheaves where $\mathcal{C}^{-i}$ is the combinatorial sheaf with
$$
\mathcal{C}^{-i}(\st{\tau}) = \{ \textrm{$i$-chains on $\st{\tau} \cap \overline{\sigma}$} \}
$$
and with the boundary maps as differentials (see \cite[\S1.12]{gm2}). From Example \ref{dual example} we see this is precisely $F\circ T (C_\sigma)$.
\end{proof}

Suppose now, in addition to $R$ being regular, Noetherian and of finite Krull dimension, that $2$ is invertible so that the constructible Witt groups form a generalised homology theory by Theorem \ref{ght}. \begin{theorem}
$F$ induces a natural transformation $H_*(-;\mathbb{L}^\cdot(R)) \to W^c_*(-)$.
\end{theorem}
\begin{proof}
We sketch the proof leaving the reader to check the details. Recall that the free symmetric L-group $H_n(K; \mathbb{L}^\cdot(R))$ is given by the cobordism classes of $n$-dimensional symmetric Poincar\'e complexes in $\comrxmod$. Using Lemma \ref{duality lemma} we see that an $n$-dimensional symmetric Poincar\'e complex maps under $F$ to a complex of sheaves in $\constr{K_S}$ equipped with a morphism to the $(-n)^\textrm{th}$ shift of its Poincar\'e--Verdier dual. The conditions in Definition \ref{spc} guarantee that this map will be a symmetric isomorphism. Thus it generates a class in $W^c_n(K_S)$. Furthermore we can check that cobordant symmetric Poincar\'e complexes give rise to the same class in $W^c_n(K_S)$.

Naturality follows from the last part of Lemma \ref{fff lemma}.
\end{proof}

\begin{corollary}
Assume that $R$ is a regular Noetherian ring of finite Krull dimension in which $2$ is invertible and further that we can resolve any finitely generated $R$-module by a finite complex of finitely generated free $R$-modules. Then for a finite simplicial complex $K$ there is a natural isomorphism $H_*(K;\mathbb{L}^\cdot(R)) \to W^c_*(K)$.
\end{corollary}
\begin{proof}
The free symmetric L-groups are a generalised homology theory. By Theorem \ref{ght} the assumption on $R$ guarantees that  the constructible Witt groups are defined and are also a generalised homology theory. Since $K$ is a finite simplicial complex it is sufficient to check that we obtain an isomorphism for a point. 

Let $\cat{T}$ be the triangulated category obtained by inverting chain equivalences in the category $\textrm{Com(R,\textrm{pt})}$ of complexes of finitely generated free $R$-modules. Walter's theorem \cite[Theorem 5.3]{walter} tells us that
$$
H_*(\textrm{pt};\mathbb{L}^\cdot(R)) \cong W^{-n}(\cat{T}).
$$
Since, by assumption, every finitely generated module has a finite resolution by finitely generated free $R$-modules the inclusion of $\cat{T}$ into the derived category of finitely generated modules, \ie into $\constr{\textrm{pt}}$, is an equivalence. Hence we also have 
$
W^{-n}(\cat{T}) \cong W^{-n}(\constr{\textrm{pt}}) \cong W_{n}^c(\textrm{pt})
$
as required.
\end{proof}

The conditions of the corollary are satisfied if, for example, $R$ is a principal ideal domain, a polynomial ring over a field or a Noetherian local ring and $2$ is invertible in $R$. In particular, they are satisfied when $R=\qq$ and this is the case we study in the next section. 
 
\section{Rational coefficients, Witt spaces and L-classes}
\label{rational theory}
We work in the PL-category --- all spaces in this section are polyhedra and all maps are piecewise linear. We discuss properties of the constructible Witt groups in the special case when the ring $R$ is the rationals $\qq$. All coefficients in homology and cohomology groups are also rational.

 In \S\ref{witt bordism} we interpret the corresponding generalised homology theory $W^c_*(-;\qq)$ geometrically as the bordism theory of Witt spaces. Then, in \S\ref{L-classes}, we show that the constructible Witt groups form the natural domain of definition for L-classes.
\subsection{Witt spaces and bordism}
\label{witt bordism}

Let $X$ be an $n$-dimensional polyhedron. $X$ is a (PL) \emph{$n$-pseudomanifold} if there is a closed subspace $Y$ with $\dim Y \leq n-2$ such that $X - Y$ is an $n$-manifold which is dense in $X$. An $n$-pseudomanifold with (collared) boundary is a pair $(X,\bdy X)$ such that $X - \bdy X$ and $\bdy X$ are pseudomanifolds, of respective dimensions $n$ and $n-1$, and a neighbourhood of $\bdy X$ in $X$ is (PL)-homeomorphic to $\bdy X \times [0,1)$.
\begin{remark}
In terms of triangulations, $X$ is an $n$-pseudomanifold if, for any triangulation, $X$ is the union of the $n$-simplices and each $(n-1)$-simplex is a face of \emph{exactly} two $n$-simplices. $(X,\bdy X)$ is an $n$-pseudomanifold with boundary if, for any triangulation, $X$ is the union of the $n$-simplices, each $(n-1)$-simplex is a face of \emph{at most} two $n$-simplices, and the set of $(n-1)$-simplices which are the face of only one $n$-simplex forms an $(n-1)$-pseudomanifold $\bdy X$. 
\end{remark}

A pseudomanifold $X$ is \emph{orientable} if, and only if, $X - Y$ is orientable. If $X$ is an oriented pseudomanifold we denote the oppositely oriented pseudomanifold by $\overline{X}$.  A compact oriented $n$-pseudomanifold carries a fundamental class $[X] \in H_n(X)$ so there is a well defined pairing
\begin{equation}
\label{coh pairing}
H^i(X) \otimes H^{n-i}(X) \to H^n(X) \stackrel{\langle - ,[X]\rangle}{\longrightarrow}  \qq.
\end{equation}
Thus we can ask, for what class of pseudomanifolds is this non-degenerate \ie does Poincar\'e duality hold? Of particular interest are those $X$ for which \emph{local} Poincar\'e duality holds (from which the global version follows), namely the rational homology manifolds. 

The intersection cohomology groups $\ih{*}{X}$ of a pseudomanifold $X$  were introduced in \cite{gm1} to study singular spaces. Their key property is that there is an analogous pairing
\begin{equation}
\ih{i}{X} \otimes \ih{n-i}{X} \to H^n(X) \stackrel{\langle -  ,[X]\rangle}{\longrightarrow}  \qq.
\end{equation}
\label{int coh pairing}
for a compact, oriented $n$-pseudomanifold $X$, which is non-degenerate for a wider class of spaces than (\ref{coh pairing}). In particular there is an interesting class of spaces for which the local (and hence the global) version of this intersection cohomology Poincar\'e duality holds; these are the Witt spaces.
\begin{definition}[Siegel {\cite[I.2]{siegel}}]
The link of a point $x$ in an $n$-pseudomanifold $X$ is an $(n-1)$-pseudomanifold. It is unique up to PL-homeomorphism and is PL-homeomorphic to the join $S^{d(x)-1}*L(x)$ where $L(x)$ is a pseudomanifold, again unique up to PL-homeomorphism, of dimension $l(x)= n-d(x)-1$. 
The pseudomanifold $X$ is a \emph{Witt space} if $\ih{l(x)/2}{L(x)} = 0$ for all $x \in X$ with $l(x)$ even. A \emph{Witt space with boundary} is a pseudomanifold with boundary $(X, \bdy X)$ where $X$ is a Witt space (in which case it follows that $\bdy X$ is also a  Witt space).
\end{definition}
\begin{remark}
If we stratify $X$, for instance by choosing a triangulation, then it is a Witt space if, and only if, for every $(2k+1)$-codimensional stratum $S$ the middle dimensional intersection cohomology $\ih{k}{L(S)}$ of the link $L(S)$ of the stratum vanishes. There is no condition on strata of even codimension. 
\end{remark}

\begin{examples}
Clearly any manifold is a Witt space. Any pseudomanifold which can be stratified with only even dimensional strata, for instance any complex projective variety, is also a Witt space.  
\end{examples}

\begin{definition}
Given a pair $(X,A)$ let $\textrm{Witt}/(X,A)$ be the category of \emph{Witt spaces over $(X,A)$} whose  objects are compact oriented Witt spaces with boundary $(W,\bdy W)$,
equipped with a map of pairs $f:(W,\bdy W) \to (X,A)$. The morphisms are commuting diagrams
$$
\xymatrix{
(W,\bdy W) \ar[dr]_{f} \ar[rr] && (W',\bdy W') \ar[dl]^{f'} \\
& (X,A).
}
$$
\end{definition}

\begin{definition}
Objects $(W,\bdy W)$ and $(W',\bdy W')$ in $\textrm{Witt}/(X,A)$ which share a common, but oppositely oriented, boundary component $V$ can be glued together to form a new Witt space
$$
(W \cup_V W' , (\bdy W - V) \sqcup (\bdy W' - V) )
$$
over $(X,A)$.
\end{definition}

\begin{definition}
Isomorphism classes of Witt spaces over $(X,A)$ form a monoid under disjoint union $\sqcup$. It is graded by dimension. The \emph{bordism group of Witt spaces} $\Omega^{\textrm{Witt}}_*(X,A)$ is the quotient of this monoid by the submonoid generated by spaces $(W,\bdy W)$ such that there exist $(W',\bdy W') \in \textrm{Witt} / (A,A)$ and $(V,\bdy V) \in \textrm{Witt} / (X,X)$ with
$$\bdy W' = \overline{\bdy W}\ \textrm{and} \ W \cup_{\bdy W} W' = \bdy V.$$
This is an Abelian group: $(W,\bdy W) \times [0,1]$ is a Witt space with boundary $$W \cup_{\bdy W} (\bdy W \times [0,1]) \cup_{\overline{\bdy W}} \overline{W}.$$ Hence the class of $(W,\bdy W) \sqcup (\overline{W},\overline{\bdy W})$ is zero in the bordism group.

\end{definition}

The bordism groups of Witt spaces form a generalised homology theory in the usual way. Siegel \cite{siegel} identifies the point groups as follows. Given a $4k$-dimensional Witt space $X$ (for $k>0$) the intersection pairing 
$$
\ih{2k}{X} \otimes \ih{2k}{X} \to  \qq
$$
defines a non-degenerate symmetric rational bilinear form $I_X$. The assignment
$$
X \mapsto [I_X] \in W(\qq)
$$
descends to an isomorphism $\Omega_{4k}^{\textrm{Witt}} \cong W(\qq)$ from the bordism group of $4k$-dimen\-sional Witt spaces for $k>0$. Apart from $\Omega_0^{\textrm{Witt}}$, all the other bordism groups vanish:
$$
\Omega_i^{\textrm{Witt}} \cong \left\{
\begin{array}{ll}
\zz & i=0\\
W(\qq) & i=4k,\ k>0\\
0 & \textrm{otherwise.}
\end{array}
\right.
$$
(The appearance of the rational Witt group $W(\qq)$ explains the name Witt spaces.) 
\begin{remark}
The structure of the rational Witt group is well-known; see, for example, \cite[Example 2.8]{balmer} or the classic \cite{mh}. It is
$$
W(\qq) \cong \zz \bigoplus_{\textrm{primes}\ p} W(\mathbb{Z}_p)
$$
where
$$
W(\mathbb{Z}_p)\cong \left\{
\begin{array}{ll}
\zz_2 & p=2\\
\zz_2 \oplus \zz_2 & p=1 \mod 4\\
\zz_4 &p=3 \mod 4.
\end{array}
\right.
$$
The inital $\zz$ corresponds to the signature of the form.
\end{remark}

Siegel's result allows us to interpret the constructible Witt groups of a polyhedron geometrically. We work relative to a fixed polyhedron $X$. Let $f:Y \to X$ be a compact $i$-dimensional pseudomanifold over $X$. In \cite{gm2} Goresky and MacPherson construct an object $\ic{Y} \in \constrc{Y}$ whose hypercohomology is the  intersection cohomology of $Y$. When $Y$ is a compact oriented Witt space they also construct a symmetric Verdier dual pairing
\begin{equation}
\label{gm pairing}
\ic{Y} \otimes \ic{Y} \to \shift^{-i}\dualiser_Y
\end{equation}
or, equivalently, a symmetric isomorphism $\ic{Y} \to \shift^{-i}\dual\ic{Y}$.

\begin{remark}
\label{gm rem}
Unfortunately, there are several indexing conventions for the intersection cohomology complex $\ic{X}$ on a pseudomanifold $X$ which differ by shifts. For us $\ic{X}$ will be an extension of the constant sheaf $\qq_{U}$ on the nonsingular part $U$ of $X$ \emph{placed in degree $0$}. This is in contrast to the convention in \cite{gm2} where $\ic{X}$ is an extension of $\shift^{\dim X}\qq_{U}$. Our convention has the advantage that the $i^{th}$ intersection cohomology group is the $i^{th}$ hypercohomology of $\ic{X}$:
$$
\ih{i}{X} = H^i(X;\ic{X})
$$
For a reduced space the above symmetric isomorphism is unique up to sign (corresponding to the two choices of orientation). We refer to this as the Goresky--MacPherson isomorphism. 
\end{remark}

\begin{definition}
If $Y \to X$ is a compact oriented $i$-dimensional Witt space over $X$ then the derived pushforward of the symmetric isomorphism $\ic{Y} \to \shift^{-i}\dual\ic{Y}$ yields a representative for a class in $W^c_i(X)$ which we denote $[Y]_W$. In particular, if $X$ is itself a Witt space, then it carries a Witt orientation $[X]_W\in W^c_{\dim X}(X)$.
\end{definition}

\begin{lemma}
Suppose $Z \to X$ is a compact oriented $(i+1)$-dimensional Witt space over $X$ with boundary $\bdy Z = Y$. Then $[Y]_W=0\in W^c_i(X)$. 
\end{lemma}
\begin{proof}
This follows from Poincar\'e--Lefschetz duality for intersection cohomology. More formally, let $\imath : Z -\bdy Z \hookrightarrow Z$ and $\jmath : \bdy Z \hookrightarrow Z$ be the inclusions. By Proposition \ref{pushforward} there is a natural symmetric morphism $\rdf{\imath_!} \to \rdf{\imath_*}$ of functors $\constr{Z - \bdy Z} \to \constrc{Z}$. As in \S\ref{htpy} we can explicitly identify the cone as a symmetric isomorphism
\begin{equation}
\label{boundary cone}
\rdf{\jmath_*}\jmath^*\rdf{\imath_*} \to \shift \dual \rdf{\jmath_*}\jmath^*\rdf{\imath_*} \dual.
\end{equation}
Since $Z- \bdy Z$ is a Witt space we have a symmetric isomorphism
\begin{equation}
\label{symm str on interior}
\ic{Z - \bdy Z} \to \shift^{-i-1} \dual \ic{Z - \bdy Z}
\end{equation}
in $\constr{Z - \bdy Z}$. Standard results show that $\rdf{\jmath_*}\jmath^*\rdf{\imath_*}\ic{Z - \bdy Z} \cong \rdf{\jmath_*}\ic{\bdy Z}$. Further the symmetric isomorphism
$$
\rdf{\jmath_*}\ic{\bdy Z} \to \shift^{-i}\dual \rdf{\jmath_*}\ic{\bdy Z}
$$
arising from (\ref{boundary cone}) and (\ref{symm str on interior}) is, by the uniqueness alluded to in Remark \ref{gm rem}, the Goresky--MacPherson isomorphism for $\bdy Z$. The result follows from the second part of Proposition \ref{functoriality prop}.
\end{proof}
\begin{corollary}
For a compact polyhedron $X$ there are natural maps
\begin{equation}
\label{geometric interp}
\Omega^{\textrm{Witt}}_i(X) \to W_i^c(X)
\end{equation}
which are isomorphisms for $i>\dim X$.
\end{corollary}
\begin{proof}
The existence of the maps follows directly from the above lemma. Both $\Omega^{\textrm{Witt}}_*$ and $W^c_*$ are generalised homology theories. It follows from \cite[Theorem 5.6]{twg2} that 
$$
W^c_i(\mathrm{pt}) \cong \left\{
\begin{array}{ll}
W(\qq) & i = 0 \mod 4 \\
0 &  i \neq 0 \mod 4.
\end{array}
\right.
$$
Siegel's result then shows that (\ref{geometric interp}) holds when $X$ is a point. A standard induction over the number of simplices in a triangulation of $X$ completes the proof. 
\end{proof}
We can rephrase the connection between constructible Witt groups and Witt bordism as follows. Note that the product of two Witt spaces is also a Witt space and the (external) product structure on Witt groups discussed in \S\ref{products} arises from the evident product
\begin{eqnarray}
\label{witt product}
\textrm{Witt}/{(X,A)} \times \textrm{Witt}/{(Y,B)} &\to & \textrm{Witt}/{(X\times Y ,A\times Y \cup X \times B)}.
\end{eqnarray}
The $4$-periodicity of the constructible Witt groups arises from taking the product with the class of $\cc\pp^2$ in $\Omega^{\textrm{Witt}}_4(\textrm{pt}) \cong W^c_4(\textrm{pt})$. It  follows from the above corollary that the sequence
$$
\Omega^{\textrm{Witt}}_{i}(X)  \to \Omega^{\textrm{Witt}}_{i+4}(X) \to \Omega^{\textrm{Witt}}_{i+8}(X) \to \ldots 
$$
arising from taking products with $\cc\pp^2$ stabilises, and that
\begin{eqnarray}
\label{bordism colimit}
W_i^c(X) \cong \textrm{colim}_{k\to\infty} W_{i+4k}^c(X) \cong  \textrm{colim}_{k\to\infty} \Omega^{\textrm{Witt}}_{i+4k}(X).
\end{eqnarray}
Thus, stably, up to Witt-equivalence, every symmetrically self-dual complex of sheaves is `of geometric origin' \ie arises as the pushforward of the intersection cohomology complex on a Witt space. Bordism invariants of Witt spaces over $X$ which are stable under product with $\cc\pp^2$ correspond to  Witt-equivalence invariants of self-dual complexes of sheaves in $\constrc{X}$.

\subsection{L-classes}
\label{L-classes}

As an example of the utility of the geometric interpretation of the constructible Witt groups as bordism groups we show how it it can be used to view L-classes as homology operations from the rational constructible Witt groups to  rational homology. It should be noted that this is not the only way to proceed; a more sophisticated approach to the definition of L-classes is taken in \cite[\S 5]{cs} (based on \cite{csw}), and this could be used to give alternative proofs of the results below. We work in the PL category; see \cite{gm1} and \cite{bcs} for analogous accounts of L-classes for Whitney stratified Witt spaces.

Mimicking the approach to defining combinatorial Pontrjagin classes for rational homology manifolds in \cite[\S20]{milstaff} we obtain the following analogue of \cite[Lemmas 20.3 and 20.4]{milstaff}. (We have stated a relative version of the result but this is an easy extension, see \cite[p242]{milstaff}.)
\begin{proposition}
Let $(W,\partial W)$ be an oriented $n$-dimensional PL Witt space with boundary, $S^i$ the standard PL $i$-sphere and $f: (W,\partial W) \to (S^i,p)$ a PL map where $p \in S^i$. Then for all $q$ not in some proper closed Pl subspace of $S^i$ the fibre $f^{-1}(q)$ is an oriented PL Witt subspace (with boundary) of $W$ of dimension $\dim W - i$ and with trivial normal bundle. Furthermore, the signature $\sigma(f^{-1}(q))$ is independent of the choice of $q$ and only depends on the homotopy class of $f$.
\end{proposition}

Hence, using simplicial approximation, we obtain a well-defined map $$[(W,\bdy W),(S^i,p)] \to \zz$$ and thence, using the fact that the cohomotopy set $[(W,\bdy W),(S^i,p)]$ is a group for $2i > n+1$ and the rationalisation 
\begin{equation}
\label{hurewicz}
[(W,\bdy W),(S^i,p)] \otimes \qq \to H^i(W,\bdy W;\qq)
\end{equation}
of the Hurewicz map is an isomorphism in this range, a map
$$
H^i(W,\bdy W;\qq) \cong [(W,\bdy W),(S^i,p)] \otimes \qq \to \qq.
$$
By the same argument as \cite[Lemma 20.4]{milstaff} this is a homomorphism and so defines a class $L_i(W,\partial W)$ in $H_i(W,\partial W;\qq)$ called the $i^{th}$ L-class. The L-classes of a smooth manifold are Poincar\'e dual to the Hirzebruch L-classes of the tangent bundle, see \cite[\S 20]{milstaff}.

It follows from this definition that $L_i(W,\partial W)$ is the unique homology class such that for any normally non-singular codimension $i$ subspace $V \subset (W - \partial W)$  with trivial normal bundle $\langle \, [V] , L_i(W,\partial W) \rangle = \sigma(V)$. Here $[V] \in H^i(W,\partial W)$ is `Poincar\'e dual' to the normally non-singular subspace $V$. An important consequence (following immediately from the geometric definition of the coboundary map) is that 
\begin{equation}
\label{bcs theorem}
\partial  L_i(W,\partial W) = L_{i-1}(\partial W) \in H_{i-1}(\partial W)
\end{equation}
cf.\ the analogue for Whitney stratified Witt spaces in {\cite[\S 2]{bcs}}.

L-classes for $2i \leq n+1$ can be defined, and shown to satisfy (\ref{bcs theorem}), by taking products of $(W, \bdy W)$ with spheres.

To compare L-classes it makes sense to work relative to a fixed base \ie to work in $\textrm{Witt}/(X,A)$. The $i^{th}$ L-class defines a functor (which we rather sloppily also refer to as the $i^{th}$ L-class but denote by $\mathcal{L}_i$ rather than $L_i$)
\begin{eqnarray*}
\textrm{Witt}/(X,A) & \to&  H_i(X,A)\\
\left( f: (W,\bdy W) \to (X,A) \right)& \mapsto & f_*L_i(W).
\end{eqnarray*}
\begin{lemma}
\label{novikov add}
\begin{enumerate}
\item For $f : (W, \bdy W) \to (X,X)$ we have $\mathcal{L}_i(\bdy W) = 0$ in $H_i(X)$. In particular the L-classes are bordism invariants of Witt spaces.
\item For $(W, \bdy W)$ and $(W', \bdy W')$ in $\textrm{Witt}/(X,A)$ we have $$\mathcal{L}_i(W \cup_V W') = \mathcal{L}_i(W) + \mathcal{L}_i(W')$$ in $H_i(X,A)$.
\end{enumerate}
\end{lemma}
\begin{proof}
\begin{enumerate}
\item There is a commutative diagram of long exact sequences
$$
\xymatrix{
\ldots \ar[r] & H_{i+1}(W,\bdy W) \ar[d]^{f_*} \ar[r]^{\bdy} & H_i(\bdy W) \ar[r]\ar[d]^{f_*} & H_i(W) \ar[r] \ar[d]^{f_*}& \ldots\\
\ldots \ar[r] & H_{i+1}(X,X) \ar[r]^{\bdy} & H_i(X) \ar[r] & H_i(X) \ar[r] & \ldots
}
$$
in which $H_{i+1}(X,X) = 0$. By (\ref{bcs theorem}) we have $L_i(\bdy W) = \bdy L_{i+1}(W)$ so that $\mathcal{L}_i(\bdy W) = f_*L_i(\bdy W) = f_* \bdy L_{i+1}(W) = \bdy f_* L_{i+1}(W) =0$. 
\item Consider the evaluation of $\mathcal{L}_i(W \cup_V W')$ on a class in $H^i(X,A)$ represented by a map to $(S^i,p)$. This is given by taking the signature of the inverse image of some $q \neq p$ under a PL map in the homotopy class of the composite 
$$
\xymatrix{
(W \cup_V W', (\bdy W - V) \sqcup (\bdy W' - V)) \ar[d] \ar@{-->}[dr]^f & \\
(X,A) \ar[r] & (S^i,p).
}
$$
The fibre $g^{-1}(q)$ will be the disjoint union of the Witt subspaces $(g|_W)^{-1}(q)$ and $(g|_{W'})^{-1}(q)$ and $g|_W$ and $g|_{W'}$ will be homotopic to $f|_W$ and $f|_{W'}$ respectively. The result now follows from the fact that the signature of a disjoint union is the sum of the signatures of the components.
\end{enumerate}
\end{proof}

\begin{corollary}
$\mathcal{L}_i$ descends to a functor $\Omega^\textrm{Witt}_*(X,A) \to H_i(X,A)$.
\end{corollary}
\begin{proof}
We need to check that $\mathcal{L}_i$ vanishes on the submonoid of null-bordant spaces. Recall that  $(W,\bdy W)$ is null-bordant if there exists $(W',\bdy W') \in \textrm{Witt}/(A,A)$ with $\bdy W' = \overline{\bdy W}$, and $(V,\bdy V) \in \textrm{Witt}/(X,X)$ with $W \cup_{\bdy W} W' = \bdy V$. Then, using Lemma \ref{novikov add},
$$
0  =  \mathcal{L}_i(\bdy V) =\mathcal{L}_i(W \cup_{\bdy W} W')=\mathcal{L}_i(W) + \mathcal{L}_i(W')=\mathcal{L}_i(W).
$$
\end{proof}

\begin{proposition}
\label{L-class properties}
The L-classes form a set of stable homology operations. In other words
\begin{enumerate}
\item (naturality) they commute with pushforwards:
$$
\xymatrix{
\Omega^\textrm{Witt}_*(X,A) \ar[r]^{\mathcal{L}_i} \ar[d]_{f_*} &  H_*(X,A)\ar[d]^{f_*} \\
\Omega^\textrm{Witt}_*(Y,B) \ar[r]^{\mathcal{L}_{i}} &  H_*(Y,B)
}
$$
where $f : (X,A) \to (Y,B)$ and,
\item (stability) they commute with boundary maps: 
$$
\xymatrix{
\Omega^\textrm{Witt}_*(X,A) \ar[r]^{\mathcal{L}_i} \ar[d]_{\bdy} &  H_*(X,A)\ar[d]^{\bdy} \\
\Omega^\textrm{Witt}_{*-1}(A) \ar[r]^{\mathcal{L}_{i-1}} &  H_{*-1}(A).
}
$$
\end{enumerate}
Furthermore they are natural with respect to products in the sense that
$$
\xymatrix{
\Omega^\textrm{Witt}_*(X,A) \otimes \Omega^\textrm{Witt}_*(Y,B) \ar[d]_{\mathcal{L} \otimes \mathcal{L}} \ar[r] & \Omega^\textrm{Witt}_*(X \times Y, A \times Y \cup X\times B) \ar[d]^{\mathcal{L}}\\
H_*(X,A) \otimes H_*(Y,B) \ar@{=}[r] & H_*(X \times Y, A \times Y \cup X\times B) 
}
$$
commutes, where $\mathcal{L} = \oplus_i \mathcal{L}_i: \Omega^\textrm{Witt}_* \to H_*$ denotes the \emph{total L-class},  and the horizontal maps are that induced by (\ref{witt product}) and the Kunneth isomorphism respectively.
\end{proposition}
\begin{proof}
Naturality with respect to maps follows easily from the geometric definition. Stability is a direct consequence of (\ref{bcs theorem}). (Note that this notion of stability is equivalent to the usual definition in terms of the corresponding reduced homology theory and suspensions.) 

To prove naturality of the L-classes with respect to products we proceed as follows. For each sphere $S^i$ choose a basepoint $p^i$ and another, distinct, point $q^i$ such that $q^i \wedge q^j = q^{i+j}$ where $\wedge$ denotes the smash product. Given maps $f^i:(W,\bdy W) \to (S^i,p^i)$ and $f^j:(W',\bdy W') \to (S^j,p^j)$, representing classes respectively in $H^i(W,\bdy W)$ and $H^j(W',\bdy W')$, the composition
$$
f^{i+j} : W \times W' \to S^i \times S^j \to S^i \wedge S^j \cong S^{i+j}.
$$
of the product $f^i \times f^j$ and the quotient represents the corresponding class  in $H^{i+j}(W\times W',\bdy (W \times W'))$ under the Kunneth isomorphism. Furthermore, since $q^i$ and $q^j$ are distinct from the basepoints, we have
$$
\left(f^{i+j}\right)^{-1}(q^{i+j}) = (f^i)^{-1}(q^i) \times (f^j)^{-1}(q^j).
$$
The signature is multiplicative and it follows that
$$
\langle [f^{i+j}] , {L}_{i+j}(W \times W')\rangle = \langle [f^i] ,{L}_{i}(W )\rangle \cdot \langle [ f^j] , {L}_{j}(W')\rangle.
$$
An appeal to the Kunneth theorem completes the proof of naturality with respect to products.
\end{proof}
We make some simple remarks about L-classes.  It is clear that 
$$
\mathcal{L}_i : \Omega^\textrm{Witt}_k(X,A) \to H_i(X,A)
$$
vanishes for $i<0$ and for $i>k$. For $W \in \textrm{Witt} / (\textrm{pt},\emptyset)$ it is easy to identify $\mathcal{L}_0(W) = \sigma(W)$. For $i$-dimensional $(W,\bdy W) \in \textrm{Witt} / (X,A)$ the $i^{th}$ L-class $\mathcal{L}_i(W)$ is the image of the fundamental class of $W$ under $H_i(W,\bdy W) \to H_i(X,A)$. Finally, since rational homology is torsion-free, the L-classes of any torsion element in $\Omega^\textrm{Witt}_*(X,A)$ must vanish. Some computations of L-classes for singular spaces can be found in \cite{bcs}.

\begin{lemma}
L-classes are stable under the map $\Omega_*^{\textrm{Witt}}(X) \to \Omega_{*+4}^{\textrm{Witt}}(X)$ induced by product with $\cc\pp^2$ \ie
$$
\xymatrix{
\Omega_*^{\textrm{Witt}}(X) \ar[rr]^{\times \cc\pp^2} \ar[dr]_{\mathcal{L}_i} &&\Omega_{*+4}^{\textrm{Witt}}(X) \ar[dl]^{\mathcal{L}_i}\\
& H_i(X;\qq)
}
$$
commutes.
\end{lemma}
\begin{proof}
This follows from the naturality of the total L-class under products and $\mathcal{L}_0(\cc\pp^2) = \sigma(\cc\pp^2) = 1$.
\end{proof}
An immediate consequence is that we can define the L-classes of elements in the constructible Witt groups via 
$$
W_i^c(X) \cong \textrm{colim}_{k\to\infty} \Omega^{\textrm{Witt}}_{i+4k}(X) \stackrel{\mathcal{L}_i}{\longrightarrow} H_i(X;\qq)
$$
using (\ref{bordism colimit}). We leave the reader to prove the analogue of Proposition \ref{L-class properties}.
%\bibliographystyle{plain}
%\bibliography{witt}

\end{document}